\documentclass[journal=jacsat,manuscript=article]{achemso}

\usepackage{color}
\usepackage{float}
\usepackage{bm}
\usepackage{amsmath}
\usepackage{amssymb}
\usepackage{graphicx}
\usepackage{setspace}
\usepackage{subscript}
\usepackage{upgreek}
\usepackage{natbib}
\usepackage{babel}
\usepackage{mathtools}

\makeatletter
\newcommand{\diag}{\mathop{\mathrm{diag}}}
\DeclareMathOperator{\tr}{tr}
\DeclareMathOperator*{\argmax}{arg\,max}

\usepackage[algo2e,ruled]{algorithm2e} 
\SetKwInOut{Input}{Input}
\SetKwInOut{Output}{Output}
\mciteErrorOnUnknownfalse
\SectionNumbersOn

\title{Combining Gaussian processes and polynomial chaos expansions for stochastic nonlinear model predictive control}

\author{Eric Bradford}
\email{ecb1992@outlook.com}
\affiliation[NTNU]{Department of Engineering Cybernetics, Norwegian University of Science and Technology, Trondheim, Norway}
\alsoaffiliation[ETH]
{Institute of Environmental Engineering, ETH Zürich, Zürich, Switzerland}

\author{Lars Imsland}
\email{lars.imsland@ntnu.no}
\affiliation[NTNU]{Department of Engineering Cybernetics, Norwegian University of Science and Technology, Trondheim, Norway}

\begin{document}
\singlespacing

\begin{abstract}
Model predictive control is an advanced control approach for multivariable systems with constraints, which is reliant on an accurate dynamic model. Most real dynamic models are however affected by uncertainties, which can lead to closed-loop performance deterioration and constraint violations. In this paper we introduce a new algorithm to explicitly consider time-invariant stochastic uncertainties in optimal control problems. The difficulty of propagating stochastic variables through nonlinear functions is dealt with by combining Gaussian processes with polynomial chaos expansions. The main novelty in this paper is to use this combination in an efficient fashion to obtain mean and variance estimates of nonlinear transformations. Using this algorithm, it is shown how to formulate both chance-constraints and a probabilistic objective for the optimal control problem. On a batch reactor case study we firstly verify the ability of the new approach to accurately approximate the probability distributions required. Secondly, a tractable stochastic nonlinear model predictive control approach is formulated with an economic objective to demonstrate the closed-loop performance of the method via Monte Carlo simulations. 
\end{abstract}

\section{Introduction}
Model predictive control (MPC) was developed in the late seventies and refers to a popular control approach that has been applied in the process industry to deal with multivariable control systems with important constraints. MPC solves at each sampling time a finite horizon optimal control problem (OCP) to determine a control action to take, by exploiting a dynamic model directly. Feedback enters this process by updating the initial state at each sampling time with the available measurements. Nonlinear MPC (NMPC) refers to MPC that utilize nonlinear dynamic models, which is particularly relevant for highly nonlinear problems operated at unsteady state, such as batch processes \citep{Maciejowski2002}. While most NMPC algorithms are based on data-driven nonlinear models, NMPC using first-principles is becoming increasingly used due to the availability of more efficient optimization algorithms \citep{Biegler2010}. Commonly MPC algorithms are applied for set-point tracking, while economic MPC employs as cost function the quantity to be maximized directly, such as profits \citep{Rawlings2009}. Many dynamic models include significant uncertainties, such as parametric deviations or unaccounted disturbances. These may negatively affect the performance of the MPC algorithm and lead to constraint violations. For economic MPC the system is often driven close to its constraints \citep{Lucia2014a}. It is therefore crucial to account for significant uncertainties in the formulation of the MPC problem. 

Assuming uncertainties to be deterministic and bounded leads to robust MPC (RMPC). Min-max MPC frameworks are among the first methods proposed and focus on minimizing the cost while satisfying constraints under worst-case realizations \citep{Scokaert1998}. Min-max MPC has also been applied to nonlinear systems, for example in \citet{Chen1997}. These methods however were found to be often unable to deal with the spread of the state trajectories or be overly conservative. Tube-based approaches were subsequently developed to address these limitations, which uses a feedback controller explicitly to ensure that the real system remains in a tube computed offline centred around the nominal solution \citep{Mayne2005}. Several nonlinear tube-based MPC algorithms have been proposed including \citet{Mayne2011,Marruedo2002,Kohler}. RMPC allows for the analysis of stability and performance of the system in the worst case, which may however have a diminishingly small chance of occurrence and hence can be overly conservative. 

An alternative to RMPC is given by stochastic MPC (SMPC). In SMPC the uncertainties are described by known probability density functions (pdf). In SMPC constraints are given by either expectation or chance constraints. Therefore, SMPC allows for a controlled rate of constraint violations, which avoids the previously mentioned problem of RMPC and leads to a trade-off between risk of constraint violation and performance \citep{Mesbah2016}. Recent reviews on SMPC can be found in \citet{Farina2016,Mesbah2016}. SMPC has predominantly been developed for linear systems. An important group of algorithms are given by probabilistic tube-based approaches as counter-part to their robust versions for additive and multiplicative noise \citep{Cannon2009,Cannon2011}. Alternatively several approaches have been suggested using affine parametrization of either state or disturbance feedback \citep{Korda2014,Oldewurtel2008,Hokayem2012}, which apply exact propagation methods for the mean and covariance. Lastly, scenario-based MPC methods have been put forward that simulate separate realization of the uncertain parameters and use Monte Carlo estimations of chance constraints and objectives \citep{Schildbach2014,Prandini2012,DelaPenad2005}. Stochastic NMPC (SNMPC) has on the other hand received significantly less attention, which may be in part explained by the difficulty of propagating stochastic variables through nonlinear transformations. An exception to this is given by the case of only discrete realizations of uncertainties for which efficient algorithms have been developed using multi-stage stochastic programming \citep{Lucia2014,Patrinos2014}. These propagate each possible scenario and ensure that none of the constraints are violated. Several procedures are used in literature to analyse stochastic uncertainties of nonlinear systems, including Monte Carlo (MC) sampling, polynomial chaos expansions (PCE), Gaussian closure, equivalent linearization, and stochastic averaging \citep{Konda2011}. 

A straight-forward approach for SNMPC is given by successive linearization, such as in \citet{Nagy2003} which uses an extended Kalman filter approach to propagate the stochastic uncertainties or as in \citet{Cannon2009} that applies the probabilistic tube-based approach on the successively linearized system. An alternative is given by applying the Unscented transformation \citep{Bradford2018}. Both linearization or Unscented transformations, while being computationally cheap, are only applicable to moderately nonlinear systems. \citet{Rostampour2015,Bradford2017a} utilize a sampling average approach to approximate chance constraints and objective. The required number of samples to obtain accurate predictions however quickly becomes prohibitive. In \citet{Sehr2017} an output feedback SNMPC approach is introduced by using the particle filter equations for both prediction of future state distributions and for updating the states from available measurements. Again the required number of samples can be prohibitive. \citet{Maciejowski2007} proposed to use Markov chain MC sampling to propagate the uncertainties and solve the optimization problem. The computational cost of the suggested method is high, since it aims to find the global optimum at each sampling instant. In the case of continuous time the Fokker-Planck partial differential equation system can be used to describe the evolution over time of the pdfs of the states, which is used in \citet{Buehler2016a} for SNMPC. A Lyapunov function is included to guarantee stochastic stability.   
Much of the work in SNMPC has been concerned with the application of PCEs. PCEs are an efficient alternative to MC sampling for approximating the probability distribution of a nonlinear transformation of stochastic variables by employing orthogonal polynomials \citep{Nagy2007}. PCEs in this context are a scenario-based SNMPC algorithm that uses least-squares estimation online for every iteration of inputs to approximate the coefficients of an orthogonal polynomial approximation, known as non-intrusive PCE \citep{Fagiano2012}. For polynomial-type systems Galerkin projection is used instead to determine the coefficients, which is called intrusive PCE \citep{Streif2014}. Chance constraints can either be given using Chebychev's inequality \citep{Mesbah2014} or applying a MC sampling approximation on the orthogonal polynomials themselves \citep{Streif2014}. The PCE based SNMPC algorithm has been extended to the case of output feedback in \citet{Bradford2019d,Bradford2019a} by combining the approach with a PCE nonlinear state estimator. The usefulness and generality of PCE can be seen for example by its use in \citet{Bavdekar2016a} to formulate a SNMPC formulation for design of experiments online to maintain the dynamic model or in \citet{Heirung2017} for discriminating between different dynamic models for fault-diagnosis. 

While PCE leads to useful SNMPC algorithms, it does have a few disadvantages: \vspace{-6pt} \begin{itemize} 
\itemsep0em  
\item{Computational complexity grows exponentially with the number of uncertain parameters.}
\item{Orthogonal polynomials of high-order are prone to unstable swings.}
\item{Time-varying disturbances are difficult to handle.}
\item{Expansion is only exact for infinitely many terms.}
\end{itemize}

In the statistics community PCEs are rarely used. Gaussian processes (GP) are employed instead for uncertainty analysis in "Bayesian calibration" \citep{Kennedy2001,OHagan2006}. A comparison of GPs to PCEs is given in \citet{OHagan2013}. Gaussian processes are stochastic processes that are used as non-parametric models, which unlike other popular regression methods such as neural networks not only provide predictions, but also prediction uncertainties. GPs have been applied in several MPC papers. In \citet{Kocijan2003} GPs were used to identify a nonlinear dynamic model from data as the prediction model in a NMPC algorithm. This methodology has been successfully applied in \citet{Likar2007} to control a gas-liquid separation plant. Further, in \citet{Kim2017} the GP models are updated efficiently online using recursive formulas and sparse approximations. Furthermore, GPs have been shown to be a powerful tool to model disturbances. In \citet{Maciejowski2013} a GP is used to correct a dynamic model online for fault-tolerant control, while in \citet{Klenske2016} a GP is employed to learn a function for an  unmodelled periodic error term. Similarly, \citet{Hewing2017} proposes to model residual model uncertainty by GPs. In \citet{Bradford2019e,Bradford2019h} a GP-based algorithm is proposed that tightens the constraints offline to maintain feasibility online. GPs have been in addition employed in multiple works to approximate the mean and covariance of nonlinear transformations for the Kalman filter equations \citep{Deisenroth2009,Pruher2017,Sarkka2015}, which bear some similarity to this paper's use of GPs. 

In this paper we propose a new method using GPs and PCEs jointly for SNMPC. In this regard we employ PCEs as mean function for the GP. The combination will be referred to as "GPPCE".  GPs are well-known to approximate the function well locally, but not as well globally. PCEs on the other hand are better suited for global function approximations, but may not have the same accuracy between data-points as GPs \citep{Owen2017}. The combination of both is therefore beneficial. Another advantage over regular PCEs apart from better local fits is that the uncertainty introduced through the sample approximation can be taken into account from the GPPCE, which is otherwise ignored. Furthermore, GPs are not prone to unstable swings and are interpolating, i.e. pass through all sample points exactly. Otherwise GPs suffer from similar drawbacks as PCEs. Combining GPs and PCEs for uncertainty quantification has been previously proposed in \citet{Schobi2015}. The main novelty in this paper is to show how to use this GPPCE to obtain cheap approximations of mean and variance using closed-form expressions derived in this paper. In addition, terms are identified that can be calculated $\textit{offline}$ to significantly reduce computational costs $\textit{online}$. For SNMPC the terms are utilised directly in the optimization problem and hence it is paramount that the mean and variance estimator are fast. Lastly, we show how the GPPCE expressions can be utilised to approximate the SNMPC problem. Using GPs for SNMPC was first introduced in \citet{Bradford2018b}. The remainder of the paper is structured as follows. In the first three sections of the paper we show how GPPCE can be formulated and used in an efficient fashion to propagate uncertainties. In Section \ref{sec:GP_PCE} GPs with a PCE mean function (GPPCE) are introduced. Thereafter, in Section \ref{sec:posterior_mean_variance} terms are derived to obtain posterior mean and variance estimates given a noisy input. Section \ref{sec:unc_prop_GPPCE} shows how these expressions can be utilised efficiently to propagate uncertainties. Next we show how GPPCE can be exploited to formulate a SNMPC algorithm. Section \ref{sec:prob_def} defines the general problem to be solved using the GPPCE SNMPC algorithm, while Section \ref{sec:GP_PCE_SNMPC} introduces the GPPCE SNMPC algorithm to accomplish this task. A challenging semi-batch reactor case study is outlined in Section \ref{sec:case_study}. Results and discussions to this case study are presented in Section \ref{sec:res_disc}. The paper is concluded in Section \ref{sec:conclusions}.    

\section{Problem formulation} \label{sec:prob_def}
We consider a general discrete-time stochastic nonlinear dynamic equation system with parametric uncertainties:
\begin{align}
    & \mathbf{x}_{t+1} = \mathbf{f}(\mathbf{x}_t,\mathbf{u}_t,\bm{\uptheta}), \quad \mathbf{x}_0 = \hat{\mathbf{x}}_0 \label{eq:f_x}
\end{align}
where $t$ is the discrete time, $\mathbf{x}_t \in \mathbb{R}^{n_{\mathbf{x}}}$ represents the states, $\mathbf{u}_t \in \mathbb{R}^{n_{\mathbf{u}}}$ denotes the control inputs, $\bm{\uptheta} \in \mathbb{R}^{n_{\bm{\uptheta}}}$ represents time-invariant parametric uncertainties, and $\mathbf{f}:\mathbb{R}^{n_{\mathbf{x}}} \times \mathbb{R}^{n_{\mathbf{u}}} \times \mathbb{R}^{n_{\bm{\uptheta}}} \rightarrow \mathbb{R}^{n_{\mathbf{x}}}$ is the nonlinear dynamic equation system. The parametric uncertainties $\bm{\uptheta}$ are assumed to follow a standard normal distribution, i.e. $\bm{\uptheta} \sim \mathcal{N}(\mathbf{0},\mathbf{I})$. Note this is not restrictive, since the uncertain model parameters can be parametrized in terms of $\bm{\uptheta}$ and this way obtain the required probability distribution, see for example \citet{Bradford2019a}. The initial condition is given by a known value $\hat{\mathbf{x}}_0$. Note the parametric uncertainty is assumed to be time-invariant, since PCE-based approaches are often computationally too expensive for time-varying uncertainties. The problem of time-varying noise has been previously addressed by employing conditional probability rules, which allow for the the time-varying noise to be addressed separately using for example linearization \citep{Paulson2019}. This approach could also be utilized for GPPCE to deal with time-varying noise. 

Given the dynamic system defined in Equation \ref{eq:f_x} we aim to minimize a finite horizon objective function:
\begin{subequations}
\begin{align}
& J(N,\hat{\mathbf{x}}_0,\bm{\uptheta},\mathbf{U}_N) = \mathbb{E}[J^d(N,\hat{\mathbf{x}}_0,\bm{\uptheta},\mathbf{U}_N)] + \omega \text{Var}[J^d(N,\hat{\mathbf{x}}_0,\bm{\uptheta},\mathbf{U}_N)] \\ & J^d(N,\hat{\mathbf{x}}_0,\bm{\uptheta},\mathbf{U}_N) = \mathcal{M}(\mathbf{x}_N) + \sum_{t=0}^{N-1} \mathcal{L}(\mathbf{x}_t,\mathbf{u}_t)  \label{eq:obj_def}
\end{align} 
\end{subequations}
where $N$ is the time horizon, $\mathcal{M}:\mathbb{R}^{n_{\mathbf{x}}}\rightarrow \mathbb{R}$ denotes the Mayer term, $\mathcal{L}:\mathbb{R}^{n_{\mathbf{x}} \times n_{\mathbf{u}}} \rightarrow \mathbb{R}$ represents the Lagrange term, and $\mathbf{U}_N = [\mathbf{u}_0,\ldots,\mathbf{u}_{N-1}] \in \mathbb{R}^{n_{\mathbf{u}} \times N}$ are the control actions that need to be determined.

The objective is taken as the expectation with a weighted variance added to it of a nonlinear function, i.e. the aim is to minimize the objective in Equation \ref{eq:obj_def} given the dynamic system in Equation \ref{eq:f_x}. The weighted variance can be exploited to penalize excessive uncertainty on the objective values. The case-in-point of this paper is the control of batch processes. The objective to be minimized generally depends on the final product at the end of the batch, which leads to a shrinking horizon NMPC formulation \citep{Nagy2003}. The objective depending on the final state is represented by the Mayer term.   

The control problem is subject to hard constraints on the control inputs expressed by the set $\mathbb{U}$. In addition, the control problem is subject to both nonlinear path chance constraints and terminal chance constraints. For batch processes common path chance constraints are given by safety limits, such as upper bounds on the adiabatic temperature or reactor pressure. Terminal chance constraints on the other hand often describe a minimum product quality to be reached. The constraints are formulated as follows:
\begin{subequations} \label{eq:control_con}
\begin{align} \label{eq:path_const_def}
& \mathbb{P}[g_j(\mathbf{x}_t,\mathbf{u}_t,\bm{\uptheta})\leq 0] \geq 1 - \epsilon && \forall (t,j) \in \{1,\ldots,N\} \times \{1,\ldots,n_g\} \\
& \mathbb{P}[g^N_j(\mathbf{x}_N,\bm{\uptheta})\leq 0] \geq 1 - \epsilon && \forall j \in \{1,\ldots,n_g^N\} \label{eq:terminal_const_def} \\
& \mathbf{u}_t \in \mathbb{U} && \forall t \in \{0,\ldots,N-1\}
\end{align} 
\end{subequations}
where $g_j:\mathbb{R}^{n_{\mathbf{x}} \times n_{\mathbf{u}} \times n_{\bm{\uptheta}}} \rightarrow \mathbb{R}$ are the path constraint functions, $g^N_j:\mathbb{R}^{n_{\mathbf{x}} \times n_{\bm{\uptheta}}} \rightarrow \mathbb{R}$ are the terminal constraint functions, and $\epsilon$ is the probability of constraint violations.  

The constraints are given as chance constraints due to the presence of the parametric uncertainties $\bm{\uptheta}$. Each constraint in Equations \ref{eq:path_const_def} and Equation \ref{eq:terminal_const_def} should be violated by at most a low probability $\epsilon$ despite the stochastic uncertainties present to maintain feasibility. In the following sections the GPPCE methodology is introduced to obtain cheap approximations of the mean and variance of both objective and constraint functions, which in turn can be exploited to estimate the probabilistic objective and chance constraints.    

\section{Gaussian processes with polynomial chaos expansion mean function} \label{sec:GP_PCE}
This section presents GPs and PCEs for our purposes and is not meant as a general introduction. Please refer to \citet{Rasmussen2006,Sacks1989,Jones1998} for general descriptions of GPs and refer to \citet{Kersaudy2015,OHagan2013,Nagy2007} for a general outline on PCEs. 

GP regression is utilized to identify an unknown function $\zeta:\mathbb{R}^{n_{\bm{\uptheta}}}\rightarrow \mathbb{R}$ from data. GPs are a generalization of the multivariate normal distribution and can be viewed as distributions over functions. These can hence be used as prior functions in a Bayesian framework. The posterior update of this prior then gives us the required function approximation. A GP is fully described by a mean function $m(\cdot)$ and a covariance function $k(\cdot,\cdot)$ as follows:
\begin{subequations}
\begin{align}
& m(\bm{\uptheta}) := \mathbb{E}_{\zeta}\left[\zeta(\bm{\uptheta})\right] \\
& k(\bm{\uptheta},\bm{\uptheta}') := \mathbb{E}_{\zeta}\left[(\zeta(\bm{\uptheta}) - m(\bm{\uptheta}))(\zeta(\bm{\uptheta}') - m(\bm{\uptheta}'))\right]
\end{align}
\end{subequations}
where $\bm{\uptheta},\bm{\uptheta}'$ are arbitrary inputs and $\mathbb{E}_{\zeta}[\cdot]$ denotes the expectation taken over the function space $\zeta(\cdot)$. The mean function can be seen as the "average" shape of the function, while the covariance function defines the covariance between any two output values at their corresponding inputs.  

The prior is specified by the mean function and the covariance function, which need to be chosen based on the \textit{prior knowledge} of the function to be inferred. The mean function can be chosen as any function, but in general should be chosen close to the function to be learnt. In this work we propose to use as mean function a linear regression term as in universal Kriging \citep{Kersaudy2015}:
\begin{equation}
m(\bm{\uptheta}) := \sum_{i=1}^{n_{\bm{\upphi}}} \beta_i \phi_i(\bm{\uptheta}) := \bm{\upbeta}^{\sf T} \bm{\upphi}(\bm{\uptheta}) \label{eq:mean_fcn}
\end{equation} 
where $\bm{\upbeta} \in \mathbb{R}^{n_{\bm{\upphi}}}$ is a vector containing $n_{\bm{\upphi}}$ trend coefficients $\beta_i$ and $\phi_i:\mathbb{R}^{n_{\bm{\uptheta}}} \rightarrow \mathbb{R}$ are a set of basis functions collected in $\bm{\upphi}(\bm{\uptheta})=[\phi_1(\bm{\uptheta}),\ldots,\phi_{n_{\bm{\upphi}}}(\bm{\uptheta})]^{\sf T}$. 
The exact choice of the mean function is motivated by the noise assumed on $\bm{\uptheta}$, which in this paper is for $\bm{\uptheta}$ to follow a standard normal distribution with zero mean and unit variance, i.e. $\bm{\uptheta} \sim \mathcal{N}(\mathbf{0},\mathbf{I})$. The mean function in our case are then given by multivariate Hermite polynomials. This selection of the mean function is motivated by their successful use as PCEs, see \citet{Kersaudy2015,OHagan2013,Nagy2007} for more information. These can be expressed as a tensor product of univariate Hermite polynomials of the components of $\bm{\uptheta}$:
\begin{equation}
    \bm{\upphi}_{\bm{\upalpha}} = \prod_{i=1}^{n_{\bm{\uptheta}}} \phi_{\alpha_i}(\theta_i) \label{eq:multivariate_polynomials}
\end{equation}
where $\phi_{\alpha_i}:\mathbb{R} \rightarrow \mathbb{R}$ are univariate polynomials of $\theta_i$ of order $\alpha_i$. The multidimensional index $\bm{\upalpha}=[\alpha_1,\ldots,\alpha_{n_{\bm{\uptheta}}}]$ is hence used to define the degree of each univariate polynomial and the total order of the multivariate polynomial $\upphi_{\bm{\upalpha}}$ is consequently given as $|\bm{\upalpha}|=\sum_{i=1}^{n_{\bm{\uptheta}}} \alpha_i$.

The univariate polynomials $\phi_{\alpha_i}$ are chosen to satisfy an orthogonality property according to the probability distribution of $\theta_i$, which in our case for standard normal distributions leads to Hermite polynomials:
\begin{equation}
\phi_{\alpha_i}(\theta_i) = (-1)^{\alpha_i} \exp\left(\frac{1}{2}\theta_i^2\right) \frac{d^{\alpha_i}}{d\theta_i^{\alpha_i}} \exp\left(-\frac{1}{2}\theta_i^2\right) \label{eq:hermite_poly} 
\end{equation} 

Keeping all polynomial terms up to a total order of $m$ leads then to the following expression for the mean function in Equation \ref{eq:mean_fcn}:
\begin{equation}
m(\bm{\uptheta}) = \sum_{0 \leq |\bm{\upalpha}| \leq m} \beta_{\bm{\upalpha}} \upphi_{\bm{\upalpha}}(\bm{\uptheta}) = \bm{\upbeta}^{\sf T} \bm{\upphi}(\bm{\uptheta}) \label{eq:trunc_exp}
\end{equation}
where $\bm{\upbeta} \in \mathbb{R}^{n_{\bm{\upphi}}}$ and $\bm{\upphi}(\bm{\uptheta}):\mathbb{R}^{n_{\bm{\uptheta}}} \rightarrow \mathbb{R}^{n_{\bm{\upphi}}}$ are vectors of the coefficients and polynomials of the truncated expansion respectively. The truncated series consists of $n_{\bm{\upphi}}=\frac{(n_{\bm{\uptheta}}+m)!}{n_{\bm{\uptheta}}!m!}$ terms. Note the number of terms grows exponentially with the input dimension of $\bm{\uptheta}$ and the truncation order of the polynomials.

For the covariance function we utilise the anisotropic squared-exponential (SE) \citep{Rasmussen2006}: 
\begin{equation}
k(\bm{\uptheta},\bm{\uptheta}') = \alpha^2 r(\bm{\uptheta},\bm{\uptheta}'), \quad r(\bm{\uptheta},\bm{\uptheta}') = \exp\left( -\frac{1}{2} (\bm{\uptheta}-\bm{\uptheta}')^{\sf T} \bm{\Lambda}^{-1} (\bm{\uptheta}-\bm{\uptheta}') \right) \label{eq:SE} 
\end{equation}
where $\bm{\Lambda} = \diag(\lambda_1^{2},\ldots,\lambda_{n_{\bm{\uptheta}}}^{2})$ is a diagonal matrix with $n_{\bm{\uptheta}}$ separate width scaling parameters $\lambda_i$ for each input dimension $i$ and $\alpha^2$ is the covariance magnitude. The SE is infinitely differentiable and therefore assumes the underlying function to be smooth. In addition, the SE covariance function is stationary, i.e. $k(\bm{\uptheta},\bm{\uptheta}') = k(\bm{\uptheta}-\bm{\uptheta}',\mathbf{0})$. Note the SE covariance function represents an unnormalized Gaussian pdf, see Equation \ref{eq:SE_normal}. This allows later in Section \ref{sec:posterior_mean_variance} for the derivation of closed-form expressions of the posterior mean and variance, which is not possible utilizing other covariance functions.

Let the hyperparameters of the GP prior defined in Equation \ref{eq:mean_fcn} and Equation \ref{eq:SE} be denoted as $\bm{\upxi}=[\beta_1,\ldots,\beta_{n_{\bm{\upphi}}},\alpha,\lambda_1,\ldots,\lambda_{n_{\bm{\uptheta}}}]^{\sf T}$. By choosing the mean function and covariance function the prior is now specified, however in general the hyperparameters are unknown. We therefore need to infer these from data. The data is given as \textit{noiseless samples} of the function $\zeta(\bm{\uptheta})$ at separate inputs. Given $n_s$ such responses, let $\bm{\Theta}=[\bm{\uptheta}_1,\ldots,\bm{\uptheta}_{n_s}]^{\sf T} \in \mathbb{R}^{n_s \times n_{\bm{\uptheta}}}$ be a vector of the input design and $\mathbf{z} = [\zeta(\bm{\uptheta}_1),\ldots\zeta(\bm{\uptheta}_{n_s})]^{\sf T} \in \mathbb{R}^{n_s}$ be a vector of the corresponding function values. Most commonly maximum likelihood estimation (MLE) is carried out to determine reasonable hyperparameter values. The log-likelihood of the observations $\mathbf{z}$ is:     
\begin{equation}
L(\bm{\upxi}) = -\frac{n_s}{2}\log(2\pi) - \frac{n_s}{2}\log(\alpha^2) - \frac{1}{2}\log(|\bm{\Sigma}_{\mathbf{z}}|) - \frac{\bm{\upnu}^{\sf T}\bm{\Sigma}_{\mathbf{z}}^{-1}\bm{\upnu}}{2\alpha^2}
\end{equation}
where $[\bm{\Sigma}_{\mathbf{z}}]_{ij} = r(\bm{\uptheta}_i,\bm{\uptheta}_j)$,  $\bm{\upnu} = \mathbf{z} - \mathbf{m}_{\mathbf{z}}$, $\mathbf{m}_{\mathbf{z}} = \bm{\Phi} \bm{\upbeta}$, and $\bm{\Phi} = [\bm{\upphi}(\bm{\uptheta}_1),\ldots,\bm{\upphi}(\bm{\uptheta}_{n_s})]^{\sf T} \in \mathbb{R}^{n_s \times n_{\bm{\upphi}}}$ is a matrix of the regression terms of the mean function in Equation \ref{eq:trunc_exp} evaluated at the inputs of the data. 

By setting the derivatives with respect to $\alpha^2$ and $\bm{\upbeta}$ to zero, the following closed-form expressions for the optimal MLE values of $\alpha^2$ and $\bm{\upbeta}$ as functions of $\bm{\Sigma}_{\mathbf{z}}$ and $\mathbf{z}$ can be determined \citep{Kersaudy2015}:
\begin{subequations}
\begin{align}
& \hat{\bm{\upbeta}} = \mathbf{a}^{\sf T} \mathbf{z}, \quad \mathbf{a} = {(\bm{\Phi}^{\sf T} \bm{\Sigma}_{\mathbf{z}}^{-1} \bm{\Phi})}^{-1}\bm{\Phi}^{\sf T} \bm{\Sigma}_{\mathbf{z}}^{-1}  \\
& \hat{\alpha}^2 = \frac{\bm{\upnu}^{\sf T} \bm{\Sigma}_{\mathbf{z}}^{-1} \bm{\upnu} }{n_s}
\end{align}
\end{subequations}

As pointed out in \citet{Bradford2018b} the evaluation of the scaling parameters $\lambda_i$ is too expensive for online implementations and will therefore be fixed in this work. This will however lead to a worse fit of the GP and hence a larger uncertainty with regard to the model fit. We show two different approaches for fixing this parameter in this paper. A simple but effective heuristic has been suggested in \citet{Jaakkola1999}, where all the width scales are fixed to the median of all pairwise euclidean distances in the data matrix $\bm{\Theta}$:
\begin{equation}
\hat{\lambda}_i = \text{median}(||\bm{\uptheta}_i - \bm{\uptheta}_j ||_2) \quad \forall i \in \{1,\ldots,n_{\bm{\uptheta}}\} \label{eq:median_heuristic}
\end{equation}

While this in general can lead to good solutions, it ignores the response values $\mathbf{z}$ and sets all $\lambda_i$ to the same value. In the GP-based Kalman filter the width scaling parameters are fixed instead using qualitative reasoning on the \textit{importance} of the inputs on the output from $\zeta(\cdot)$ \citep{Pruher2017}. A small $\lambda_i$ corresponds to an important input dimension $\theta_i$ to the value of $\zeta(\cdot)$, while a large value conversely indicates less significance of this input dimension. Often in these applications it is simple to generate several \textit{representative} datasets $\mathbf{z}$ \textit{offline}. 

For example, imagine we require the $\bm{\uplambda}$ values for a nonlinear function $g(\mathbf{x}_N,\mathbf{u}_N,\bm{\uptheta})$ at time $N$, such as the nonlinear objective or terminal constraints in Equation $\ref{eq:obj_def}$ and Equation \ref{eq:control_con} respectively. We then generate a possible control trajectory $\mathbf{U}=[\mathbf{u}_0,\ldots,\mathbf{u}_{N-1}]$ within the control bounds and simulate the system in Equation \ref{eq:f_x} for each value of $\bm{\uptheta}$ in the data matrix $\bm{\Theta}$. This then in turn gives us a dataset $\mathbf{z}=[g(\mathbf{x}_N,\mathbf{u}_N,\bm{\uptheta}_1),\ldots,g(\mathbf{x}_N,\mathbf{u}_N,\bm{\uptheta}_{n_s})]^{\sf T}$. From these datasets we can obtain optimal scaling values as follows \citep{Forrester2009}:
\begin{align} \label{eq:optimal_scaling}
    \hat{\bm{\uplambda}} = \argmax \left[ -\frac{n_s}{2} \log(\hat{\alpha}^2) - \frac{1}{2}\log(|\bm{\Sigma}_{\mathbf{z}}|) \right] \
 \end{align}
From these different values we can then choose values for $\hat{\bm{\uplambda}}$ that account for the \textit{importance} of the different inputs. Let the corresponding scaling matrix be given by $\hat{\bm{\Lambda}} = \diag(\hat{\lambda}_1,\ldots,\hat{\lambda}_{n_{\bm{\uptheta}}})$.

Once the hyperparameters are fixed the posterior GP is utilised to obtain predictions and corresponding uncertainty values. The posterior distribution is given by the prior distribution taking the observations $\mathbf{z}$ into account. Due to the GP prior assumptions, the observations follow a multivariate Gaussian distribution. Similarly, the value of the latent function at an arbitrary input $\bm{\uptheta}$ also follows a Gaussian distribution. The conditional distribution of $\zeta(\bm{\uptheta})$ given the observations $\mathbf{z}$ can be stated as \citep{Schobi2015}:
\begin{subequations} \label{eq:GP_posterior}
\begin{align} 
& \zeta(\bm{\uptheta})|\mathbf{z} \sim \mathcal{N} \left(m_{\zeta}(\bm{\uptheta})|\mathbf{z},\sigma_{\zeta}^2(\bm{\uptheta})|\mathbf{z} \right) \\ \label{eq:GP_posteriormean}
& m_{\zeta}(\bm{\uptheta})|\mathbf{z} = m(\bm{\uptheta}) + \mathbf{r}^{\sf T}_{\zeta,\mathbf{z}}(\bm{\uptheta}) \bm{\Sigma}_{\mathbf{z}}^{-1} (\mathbf{z} - \mathbf{m}_{\mathbf{z}})    \\ \label{eq:GP_posteriorvar}
& \sigma_{\zeta}^2(\bm{\uptheta})|\mathbf{z} = \hat{\alpha}^2\left(1 -  \bm{\kappa}^{\sf T}_{\zeta,\mathbf{z}}(\bm{\uptheta}) \mathbf{K}^{-1} \bm{\kappa}_{\zeta,\mathbf{z}}(\bm{\uptheta})\right) 
\end{align}
\end{subequations}
where $m_{\zeta}(\bm{\uptheta})|\mathbf{z}$ and $\sigma_{\zeta}^2(\bm{\uptheta})|\mathbf{z}$ are the mean and variance function of $\zeta(\bm{\uptheta})|\mathbf{z}$ at an arbitrary input $\bm{\uptheta}$ given the observations $\mathbf{z}$, $\mathbf{r}_{\zeta,\mathbf{z}}(\bm{\uptheta})=[r(\bm{\uptheta},\bm{\uptheta}_1),\ldots,r(\bm{\uptheta},\bm{\uptheta}_{n_{\mathbf{z}}})]^{\sf T}$, and $\bm{\kappa}_{\zeta,\mathbf{z}}(\bm{\uptheta}) = [\bm{\upphi}^{\sf T}(\bm{\uptheta}),\mathbf{r}^{\sf T}_{\zeta,\mathbf{z}}(\bm{\uptheta})]^{\sf T}$, and $\mathbf{K} =   \begin{bmatrix}
    \mathbf{0} & \bm{\Phi}^{\sf T}  \\
    \bm{\Phi} & \bm{\Sigma}_{\mathbf{z}} 
  \end{bmatrix}$. The mean $m_{\zeta}(\bm{\uptheta})|\mathbf{z}$ can be seen as the best-estimate of $\zeta(\bm{\uptheta})|\mathbf{z}$, while the variance $\sigma_{\zeta}^2(\bm{\uptheta})|\mathbf{z}$ can be viewed as a measure of uncertainty of this prediction.  

\begin{figure}[H]
\centering
\includegraphics[width=0.95\textwidth]{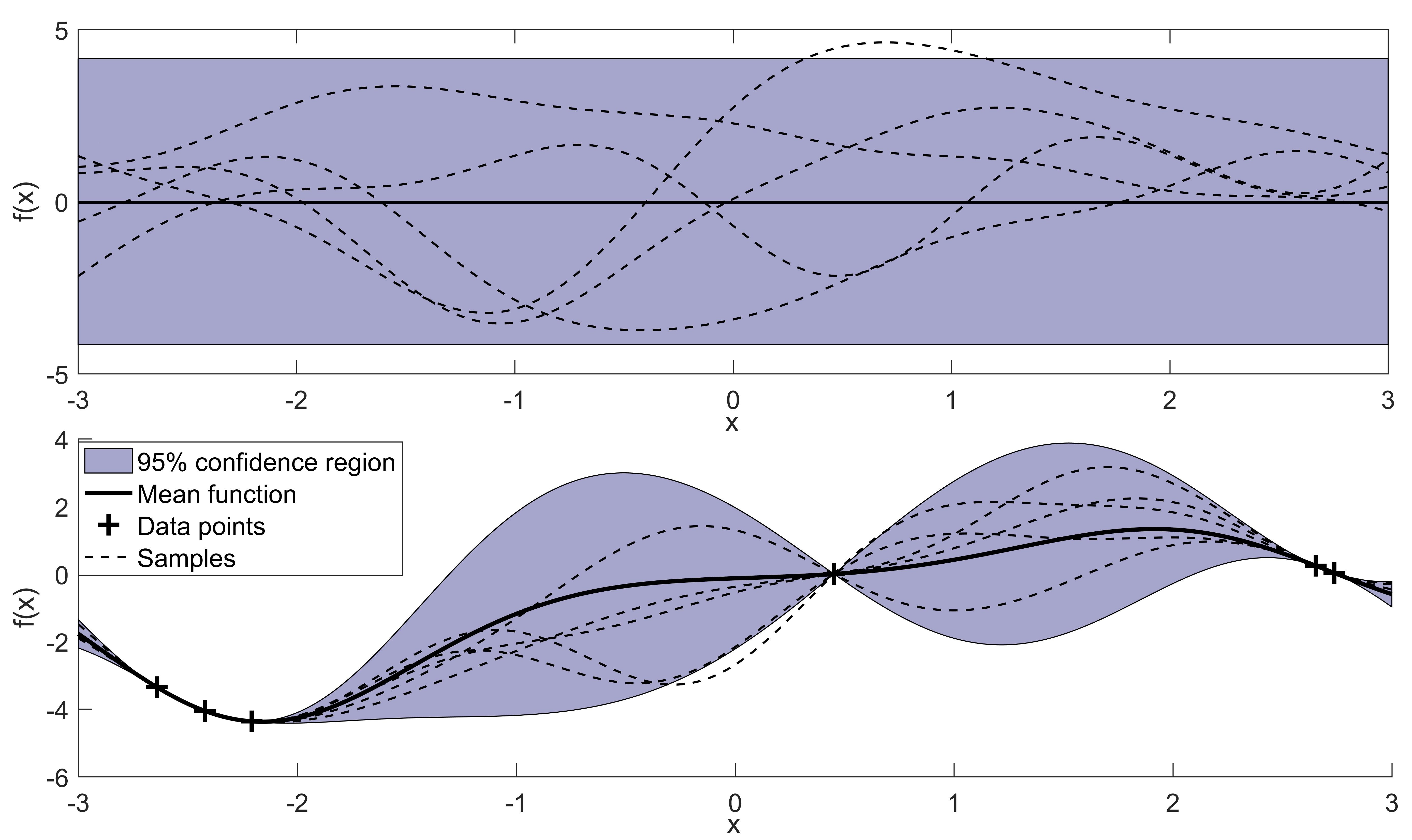}
\caption{Illustration of a GP prior is shown on the top, while the corresponding GP posterior is shown below, updated with 8 observations of a one dimensional function. The prior has a mean function of $0$ and a SE kernel given by Equation \ref{eq:SE}. It can be observed that for the GP posterior points close to the data have low uncertainty, while data far away from the observations have significantly higher uncertainty.} 
\label{fig:GP:diagram}
\end{figure}

An example of a GP prior and posterior is shown in Figure \ref{fig:GP:diagram}. Firstly, it can be seen that the posterior has significantly lower uncertainty than the prior due to the data available, especially close to the data points. Secondly, it can be seen from the samples that the SE covariance function yields smooth functions.      

\section{Posterior mean and variance estimates from GPPCE} \label{sec:posterior_mean_variance}
So far we have assumed that the input $\bm{\uptheta}$ is deterministic, often however the input $\bm{\uptheta}$ is given by a probability distribution. The aim of this section is to use the GP posterior introduced in Section \ref{sec:GP_PCE} to estimate the mean and variance of $\zeta(\bm{\uptheta})|\mathbf{z}$ given that $\bm{\uptheta}$ follows a standard normal distribution. In particular, the case of Gaussian distributed inputs has been addressed extensively due to its importance when using GP state space models for multi-step ahead predictions \citep{Girard2003}. The GPPCE surrogate $\zeta(\bm{\uptheta})|\mathbf{z}$ approximates the function $\zeta(\bm{\uptheta})$ given the observations $\mathbf{z}$, and hence its posterior mean and variance are estimates of the mean and variance of $\zeta(\bm{\uptheta})$.   

It is possible to give equations for the exact mean and variance of the GPPCE surrogate $\zeta(\bm{\uptheta})|\mathbf{z}$ for certain choices of mean and covariance function, which were made in this work. The law of iterated expectations can be used to find the exact posterior mean and variance of the GPPCE surrogate $\zeta(\bm{\uptheta})|\mathbf{z}$ as follows \citep{Quinonero-Candela}:   
\begin{subequations}
\begin{align}\label{eq:Expectation_noisymean}
& \mathbb{E}_{\bm{\uptheta}}[\zeta(\bm{\uptheta})|\mathbf{z}] = \mathbb{E}_{\bm{\uptheta}}[\mathbb{E}_{\zeta}[\zeta(\bm{\uptheta})|\mathbf{z},\bm{\uptheta}]] = \mathbb{E}_{\bm{\uptheta}}[m_{\zeta}(\bm{\uptheta})|\mathbf{z}]  \\
& \mathbb{V}_{\bm{\uptheta}}[\zeta(\bm{\uptheta})|\mathbf{z}] =  \mathbb{E}_{\bm{\uptheta}}[\mathbb{V}_{\zeta}[\zeta(\bm{\uptheta})|\mathbf{z},\bm{\uptheta}]] + \mathbb{V}_{\bm{\uptheta}}[\mathbb{E}_{\zeta}[\zeta(\bm{\uptheta})|\mathbf{z},\bm{\uptheta}]] =  \mathbb{E}_{\bm{\uptheta}}[\sigma_{\zeta}^2(\bm{\uptheta})|\mathbf{z}] + \mathbb{V}_{\bm{\uptheta}}[m_{\zeta}(\bm{\uptheta})|\mathbf{z}]
\end{align}
\end{subequations}
where $m_{\zeta}(\bm{\uptheta})|\mathbf{z}$ and $\sigma_{\zeta}^2(\bm{\uptheta})|\mathbf{z}$ are given in Equation \ref{eq:GP_posteriormean} and Equation \ref{eq:GP_posteriorvar} respectively. Note this marks a major advantage of GPPCE modelling, since this leads to analytically tractable expressions for the expectation and variance with respect to the function $\zeta(\cdot)|\mathbf{z}$ given the observations $\mathbf{z}$.  

To evaluate the posterior mean and variance we require expressions for the terms: $\mathbb{E}_{\bm{\uptheta}}[m_{\zeta}(\bm{\uptheta})|\mathbf{z}]$, $\mathbb{E}_{\bm{\uptheta}}[\sigma_{\zeta}^2(\bm{\uptheta})|\mathbf{z}]$, and $\mathbb{V}_{\bm{\uptheta}}[m_{\zeta}(\bm{\uptheta})|\mathbf{z}]$. By substituting the definitions of $m_{\zeta}(\bm{\uptheta})|\mathbf{z}$ and $\sigma_{\zeta}^2(\bm{\uptheta})|\mathbf{z}$ in Equation \ref{eq:GP_posterior}, we arrive at:
\begin{subequations} \label{eq:posterior_GP_mean_variance}
\begin{align}
\mathbb{E}_{\bm{\uptheta}}[m_{\zeta}(\bm{\uptheta})|\mathbf{z}]      & = \mathbb{E}_{\bm{\uptheta}}[m(\bm{\uptheta})] + \mathbb{E}_{\bm{\uptheta}}\left[\mathbf{r}^{\sf T}_{\zeta,\mathbf{z}}(\bm{\uptheta})\right]\bm{\Sigma}_{\mathbf{z}}^{-1} \bm{\upnu} \\
\mathbb{E}_{\bm{\uptheta}}[\sigma_{\zeta}^2(\bm{\uptheta})|\mathbf{z}] &= \hat{\alpha}^2\left(1 - \tr\left(  \mathbb{E}_{\bm{\uptheta}}\left[\bm{\kappa}_{\zeta,\mathbf{z}}(\bm{\uptheta}) \bm{\kappa}_{\zeta,\mathbf{z}}^{\sf T}(\bm{\uptheta}) \right] \mathbf{K}^{-1} \right) \right) \\
\mathbb{V}_{\bm{\uptheta}}[m_{\zeta}(\bm{\uptheta})|\mathbf{z}] &=   \mathbb{E}_{\bm{\uptheta}}\left[m(\bm{\uptheta})^2\right] + 2\mathbb{E}_{\bm{\uptheta}}\left[m(\bm{\uptheta}) \mathbf{r}^{\sf T}_{\zeta,\mathbf{z}}(\bm{\uptheta})\right] \bm{\Sigma}_{\mathbf{z}}^{-1} \bm{\upnu} + \\ & \bm{\upnu}^{\sf T} \bm{\Sigma}_{\mathbf{z}}^{-1} \mathbb{E}_{\bm{\uptheta}}\left[\mathbf{r}_{\zeta,\mathbf{z}}(\bm{\uptheta}) \mathbf{r}^{\sf T}_{\zeta,\mathbf{z}}(\bm{\uptheta})\right] \bm{\Sigma}_{\mathbf{z}}^{-1} \bm{\upnu} - \left(\mathbb{E}_{\bm{\uptheta}}[m_{\zeta}(\bm{\uptheta})]\right)^2  \nonumber
\end{align}       
\end{subequations}

Note these expressions are given by a series of expectations and variances on the covariance function and mean function. The idea here is to choose these such that the integrals given above can be evaluated exactly. The expressions for the expectations can be found in Appendix A. Substituting these values into Equation \ref{eq:posterior_GP_mean_variance}:
\begin{subequations}
\begin{align}
    \mathbb{E}_{\bm{\uptheta}}[m_{\zeta}(\bm{\uptheta})|\mathbf{z}] &= \mu_m + \bm{\upmu}^{\sf T}_{\mathbf{r}_{\zeta,\mathbf{z}}} \bm{\Sigma}_{\mathbf{z}}^{-1} \bm{\upnu} \\
    \mathbb{E}_{\bm{\uptheta}}[\sigma_{\zeta}^2(\bm{\uptheta})|\mathbf{z}] &= \hat{\alpha}^2\left(1 - \tr\left( \mathbf{M}_{\bm{\kappa}_{\zeta,\mathbf{z}} \bm{\kappa}^{\sf T}_{\zeta,\mathbf{z}}} \mathbf{K}^{-1} \right) \right) \\
    \mathbb{V}_{\bm{\uptheta}}[m_{\zeta}(\bm{\uptheta})|\mathbf{z}] &= 
    \mu_{m^2} + 2\bm{\upmu}^{\sf T}_{m \mathbf{r}_{\zeta,\mathbf{z}}} \bm{\Sigma}_{\mathbf{z}}^{-1} \bm{\upnu} + \bm{\upnu}^{\sf T} \bm{\Sigma}_{\mathbf{z}}^{-1} \mathbf{M}_{\mathbf{r}_{\zeta,\mathbf{z}}\mathbf{r}^{\sf T}_{\zeta,\mathbf{z}}} \bm{\Sigma}_{\mathbf{z}}^{-1} \bm{\upnu} -\left(\mathbb{E}_{\bm{\uptheta}}[m_{\zeta}(\bm{\uptheta})]\right)^2
\end{align}
\end{subequations}

\section{Uncertainty propagation using GPPCE} \label{sec:unc_prop_GPPCE}
In this section we outline how the GPPCE methodology outlined can be used to efficiently propagate uncertainties through nonlinear functions. Let an arbitrary nonlinear function $\zeta(\cdot)$ be given by:
\begin{equation}
    z = \zeta(\bm{\uptheta}) \label{eq:pdef}
\end{equation}
where $\bm{\uptheta} \sim \mathcal{N}(\bm{\uptheta};\mathbf{0},\mathbf{I})$ follows a standard normal distribution. 

The aim of this section is to estimate the mean and variance of $z$ using GPs as introduced in Section \ref{sec:GP_PCE}. The estimate should be as computationally cheap as possible, since it is used \textit{online}. Therefore, the section is divided into two parts: $"\textit{Offline computation}"$ and $"\textit{Online computation}"$. $"\textit{Offline computation}"$ outlines terms that do not directly depend on the response values $z$ and can hence be determined offline based on the sample design of $\bm{\uptheta}$ alone to save significant computational time, while $"\textit{Online computation}"$ shows how to obtain the posterior mean and variance estimates given the pre-computed terms. 

\subsection{Offline computation}
First we need to decide on the number of samples $n_s$ for the approximation. Thereafter, a sample design denoted by $\bm{\Theta} = [\bm{\uptheta}_1,\ldots,\bm{\uptheta}_{n_s}]^{\sf T} \in \mathbb{R}^{n_{\mathbf{s}} \times n_{\bm{\uptheta}}}$ needs to be chosen. This sample design should lead to a reasonable function approximation of $\zeta(\cdot)$ in regions that have significant probability densities. Regions of diminishing probability densities do not require good function approximations, since they do not contribute to the expectation values. This is accomplished by ensuring that the sample design is generated according to a standard normal distribution. The most obvious, but arguably worst approach uses crude MC to obtain these samples, which however may lead to poor convergence.  Alternatively, so-called space-filling designs can be used to generate the necessary points in a unit hypercube $[0,1]^{n_{\bm{\uptheta}}}$. These sample designs can then be converted to follow a standard normal distribution by using the probit function, see for example \citet{Stein1987}. Popular space-filling sampling designs include Latin hypercube designs \citep{Stein1987}, or Quasi MC designs, such as Sobol \citep{Sobol2001}. 

Given this sample design $\bm{\Theta}$ we next set the hyperparameters $\lambda_i$ according to Section $\ref{sec:GP_PCE}$, which gives us $\hat{\bm{\Lambda}}$. These are treated differently than the other hyperparameters, since they do not have a closed-form solution and hence cannot be evaluated \textit{online} without incurring large computational times.

Next there are several terms that only depend on the sample design $\bm{\Theta}$ and the hyperparameters $\hat{\lambda}_i$, which can hence be pre-computed. These are as follows:
\begin{subequations}
\begin{align}
    & \mathbf{a} = {(\bm{\Phi}^{\sf T} \bm{\Sigma}_{\mathbf{z}}^{-1} \bm{\Phi})}^{-1}\bm{\Phi}^{\sf T} \bm{\Sigma}_{\mathbf{z}}^{-1} \\
    & \mathbf{b}  = \bm{\upmu}^{\sf T}_{\mathbf{r}_{\zeta,\mathbf{z}}} \bm{\Sigma}_{\mathbf{z}}^{-1}    \\
    &  c = \tr\left(\mathbf{M}_{\bm{\kappa}_{\zeta,\mathbf{z}} \bm{\kappa}^{\sf T}_{\zeta,\mathbf{z}}} \mathbf{K}^{-1} \right) \\
    & \mathbf{d} = \bm{\upmu}_{m \mathbf{r}_{\zeta,\mathbf{z}}}^{\sf T} \bm{\Sigma}_{\mathbf{z}}^{-1} \\
    & \mathbf{E} = \frac{\bm{\Sigma}_{\mathbf{z}}^{-1}}{n_s}            \\
    & \mathbf{F} = \bm{\Sigma}_{\mathbf{z}}^{-1} \mathbf{M}_{\mathbf{r}_{\zeta,\mathbf{z}}\mathbf{r}^{\sf T}_{\zeta,\mathbf{z}}} \bm{\Sigma}_{\mathbf{z}}^{-1} 
\end{align} \label{eq:pre_computations}
\end{subequations}
where the required parameters can be evaluated using their definitions in Section \ref{sec:GP_PCE} and Appendix A for the expectations represented by $\bm{\upmu}_{\mathbf{r}_{\zeta,\mathbf{z}}}$, $\mathbf{M}_{\bm{\kappa}_{\zeta,\mathbf{z}} \bm{\kappa}^{\sf T}_{\zeta,\mathbf{z}}}$, $\bm{\upmu}_{m \mathbf{r}_{\zeta,\mathbf{z}}}$, and $\mathbf{M}_{\mathbf{r}_{\zeta,\mathbf{z}}\mathbf{r}^{\sf T}_{\zeta,\mathbf{z}}}$.

\subsection{Online computation}
Given these pre-computed values in Equation \ref{eq:pre_computations} we then have efficient formulas to estimate mean and variance using the posterior GP. For this we need to evaluate the function $\zeta()$ at the points defined in the sample design $\bm{\Theta}$ according to Equation \ref{eq:pdef}, which gives us the response vector $\mathbf{z}=[\zeta(\bm{\uptheta}_1),\ldots,\zeta(\bm{\uptheta}_{n_{\bm{\uptheta}}})]^{\sf T}$. The posterior GP then represents a fitted model using the data in $\mathbf{z}$, for which we have efficient formulas to obtain the mean and variance. Note that this is the first time we actually use the function $\zeta()$ and hence can be repeated for different functions $\zeta()$ without incurring too high computational costs. Based on the pre-computed values we obtain the following estimates for the mean and variance of the nonlinear transformation defined in Equation \ref{eq:pdef}:  
\begin{subequations} \label{eq:estimte_eqs}
\begin{align}
    & \hat{\bm{\upbeta}} = \mathbf{a}^{\sf T} \mathbf{z} \\
    & \bm{\upnu} = \mathbf{z} - \hat{\bm{\upbeta}} \bm{\Phi}, \quad \bm{\Phi} = [\bm{\upphi}(\bm{\uptheta}_1),\ldots,\bm{\upphi}(\bm{\uptheta}_{n_s})]^{\sf T}  \\
    & \hat{\alpha}^2 = \bm{\upnu}^{\sf T} \mathbf{E} \bm{\upnu} \\
    & \mathbb{E}_{\bm{\uptheta}}[\zeta(\bm{\uptheta})] \approx \mu^z_{GP}(\cdot;\bm{\tau}) = \mu_m +  \mathbf{b} \bm{\upnu} \\
    & \mathbb{V}_{\bm{\uptheta}}[\zeta(\bm{\uptheta})] \approx \sigma^2_{m_{\zeta}}(\mathbf{z};\bm{\uptau}) = 
     \mathbb{V}_{\bm{\uptheta}}[m_{\zeta}(\bm{\uptheta})|\mathbf{z}] = \mu_{m^2} + 2\mathbf{d}\bm{\upnu} + \bm{\upnu}^{\sf T} \mathbf{F} \bm{\upnu} - (\mu_m +  \mathbf{b} \bm{\upnu})^2  \\
    &  \sigma^2_{GP}(\mathbf{z};\bm{\uptau}) = \mathbb{V}_{\bm{\uptheta}}[\zeta(\bm{\uptheta})|\mathbf{z}] = \sigma^2_{m_{\zeta}}(\mathbf{z};\bm{\uptau}) + \hat{\alpha}^2(1-c) 
\end{align}
\end{subequations} 
where $\mu^z_{GP}(\cdot;\bm{\tau})$ is the mean estimate for the GPPCE, $\sigma^2_{m_{\zeta}}(\cdot;\bm{\uptau})$ is the variance estimate of the mean function from the GPPCE and can be seen as the best-estimate of the \textit{true} variance from the GPPCE, while $\sigma^2_{GP}(\cdot;\bm{\uptau})$ is the variance of the GP accounting for the uncertainty due to using only a $\textit{finite sample}$  approximation. The variance estimate $\sigma^2_{GP}(\cdot;\bm{\uptau})$ therefore has a larger variance by adding the term $\hat{\alpha}^2(1-c)$. The variable $\bm{\uptau} = \{\bm{\Theta},\hat{\bm{\uplambda}},n_{PCE}\}$ summarizes the different choices made, that define the mean and variance estimate outside of the direct dependency on the training data $\mathbf{z}$. Firstly, the sample design $\bm{\Theta}$ and the truncation order of the PCE $n_{PCE}$ should be chosen. Thereafter, the hyperparameters $\hat{\bm{\uplambda}}$ need to be determined either using heuristics or available data as shown in Section \ref{sec:GP_PCE}. From this the terms in Equation \ref{eq:pre_computations} can be pre-computed and used for different values of $\mathbf{z}$ to obtain estimates of the mean and variance.   

\section{GPPCE stochastic nonlinear model predictive control} \label{sec:GP_PCE_SNMPC}
In this section we introduce the GPPCE based SNMPC algorithm to solve the problem outlined in Section \ref{sec:prob_def} employing the dynamic equation system in Equation \ref{eq:f_x}. Assume we are at time $t$ and we have a full state measurement $\mathbf{x}_t$ of the current state. The GPPCE equations are utilised in the optimization algorithm to obtain accurate estimates of the mean and variances of both objective and constraint functions to approximate the probabilistic objective and chance constraints.

\subsection{Uncertainty propagation} \label{sec:uncertainty_prop}
For the uncertainty propagation we apply the results outlined in Section \ref{sec:unc_prop_GPPCE}. Our aim is to approximate the mean and variance of the objective and chance constraints to formulate the GPPCE SNMPC optimization problem. For this we first create a sample design $\bm{\Theta}=[\bm{\uptheta}_1,\ldots,\bm{\uptheta}_{n_s}]^{\sf T}$. Each realization of $\bm{\uptheta}$ then represents its own nonlinear dynamic equation system:
\begin{align} \label{eq:sample_f}
    & \Tilde{\mathbf{x}}_{k+1}^{(s)} = \mathbf{f}(\Tilde{\mathbf{x}}_k^{(s)},\Tilde{\mathbf{u}}_k,\bm{\uptheta}_s), \quad \Tilde{\mathbf{x}}^{(s)}_0 = \hat{\mathbf{x}}_0 && \forall s \in \{1,\ldots,n_s\} 
\end{align}
where $\Tilde{\mathbf{x}}^{(s)}_k$ and $\Tilde{\mathbf{u}}^{(s)}_k$ denotes the states and control inputs for the sample $s$.    

Using Equation \ref{eq:sample_f} we then have separate state values for each $\bm{\uptheta}$, for which we obtain different values for the constraint functions and objective: 
\begin{subequations} \label{eq:estimate_data}
\begin{align}
    & \mathbf{z}_{J^d} = [J^d(N,\hat{\mathbf{x}}_0,\bm{\uptheta}_1,\mathbf{U}_N),\ldots,J^d(N,\hat{\mathbf{x}}_0,\bm{\uptheta}_{n_s},\mathbf{U}_N)]^{\sf T} \\ 
    & \mathbf{z}_{g_j}^{(k)} = [g_j(\Tilde{\mathbf{x}}_k^{(1)},\Tilde{\mathbf{u}}_k,\bm{\uptheta}_1),\ldots,g_j(\Tilde{\mathbf{x}}_k^{(n_s)},\Tilde{\mathbf{u}}_k,\bm{\uptheta}_{n_s})]^{\sf T} && \forall (k,j) \in \{1,\ldots,N\} \times \{1,\ldots,n_g\}     \\
    & \mathbf{z}_{g_j^N} = [g^N_j(\Tilde{\mathbf{x}}_N^{(1)},\bm{\uptheta}_1),\ldots,g^N_j(\Tilde{\mathbf{x}}_N^{(n_s)},\bm{\uptheta}_{n_s})]^{\sf T} && \forall j \in \{1,\ldots,n_g^N\}
\end{align}
\end{subequations}

Using the sample design $\bm{\Theta}$ we define the mean and variance estimate GPPCE functions by determining the hyperparameters $\hat{\bm{\uplambda}}$ and setting the truncation order for the PCE $n_{PCE}$, which defines $\bm{\uptau}=\{\bm{\Theta},\hat{\bm{\uplambda}},n_{PCE}\}$. The hyperparameters $\hat{\bm{\uplambda}}$ in general will be set to different values for the constraint and objective functions, since the inputs have varied importance. The posterior mean and variance estimates of the constraints and objective are then given by the mean and variance function as defined in Equation \ref{eq:estimte_eqs}. 

\begin{figure}[H]
\centering
\includegraphics[width=0.95\textwidth]{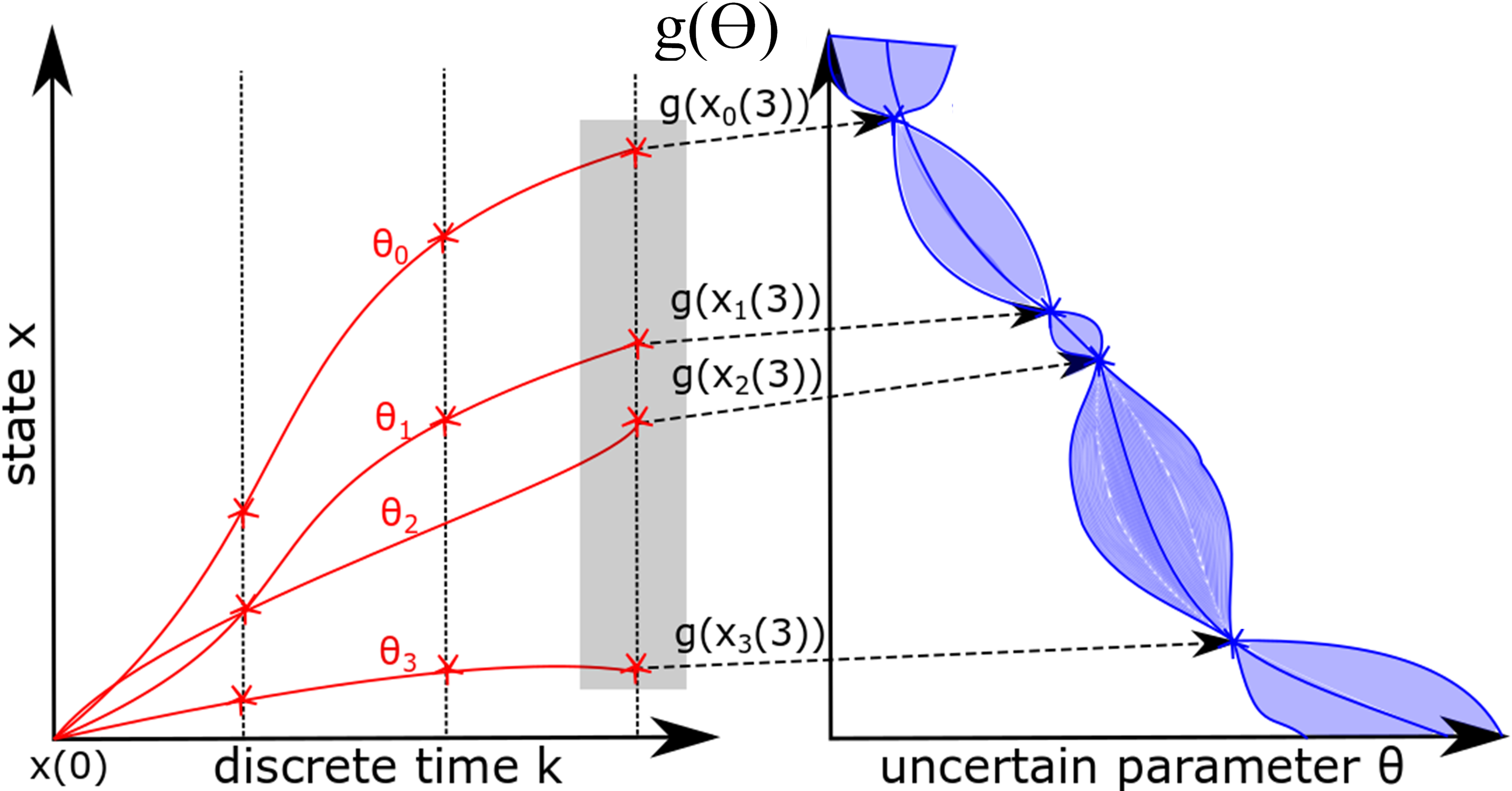}
\caption{Illustration of the mean and variance estimate from the GPPCE algorithm. The red trajectories on the left-hand side graph represent different realizations of $\bm{\uptheta}$, which each lead to different state values. These are then transformed using a constraint function $g(\cdot)$ to obtain a series of values, which are used as data to build a GP model as shown on the right-hand side graph.} 
\label{fig:GP_SNMPC_diagram}
\end{figure}

The principle of the GPPCE estimates is shown in Figure \ref{fig:GP_SNMPC_diagram}. Each sample of $\bm{\uptheta}$ corresponds to a separate trajectory according to Equation \ref{eq:sample_f}, which is shown by the red lines. Each of these trajectories then leads to distinct values of the state $\mathbf{x}$ at each discrete-time $k$ shown by the red markers. These in turn are then transformed using the objective and constraint definitions to obtain the "data" required for the mean and variance estimates, which leads to the data vectors shown in Equation \ref{eq:estimate_data}. This is highlighted by the arrows for a particular constraint. Thereafter, GP regression is applied leading to the blue line. The closed-form expressions from Equation \ref{eq:estimte_eqs} are thereafter applied, which return the exact mean and variance of the GP surrogate. Note that for each iteration of the control vector $\mathbf{U}_N$ this procedure needs to be repeated, i.e. for each step in the optimization algorithm. GPs are probabilistic models and hence also include a confidence region, which can be accounted for in the variance estimate.   

\subsection{Chance constraint approximation}
In Section \ref{sec:uncertainty_prop} we show how to obtain mean and variance estimates employing GPPCE for the constraint and objective function. These can then be utilized directly to approximate the objective in Equation \ref{eq:obj_def}. The GPPCE is a nonlinear transformation, and its output probability distribution is therefore not Gaussian. The mean and variance estimates therefore correspond to an unknown probability distribution, which necessitates the chance constraints to be reformulated using the mean and variance alone. 

The chance constraints in Equation \ref{eq:path_const_def} and Equation $\ref{eq:terminal_const_def}$ are robustly reformulated using Chebychev's inequality, since the exact evaluation of the chance constraint in Equation \ref{eq:chance_constraint_cheb} is notoriously difficult due to the integral definition of the probability function. Assume we have a chance constraints on an arbitrary random variable $\gamma$:
\begin{equation} \label{eq:chance_constraint_cheb}
    \mathbb{P}\{\gamma \leq 0\} \geq 1-\epsilon
\end{equation}

Chebychev's inequality can then be used as follows to robustly reformulate the probability constraint in Equation  \ref{eq:chance_constraint_cheb} \citep{Mesbah2014}:  
\begin{equation}
    \mu_{\gamma} + \kappa_{\epsilon} \sqrt{\sigma^2_{\gamma}} \leq 0, \quad \kappa_{\epsilon} = \sqrt{\frac{1-\epsilon}{\epsilon}}
\end{equation}
where $\mu_{\gamma}$ and $\sigma^2_{\gamma}$ are the mean and variance of $\gamma$ respectively. The robust reformulation now only requires the mean and standard deviation of $\gamma$. 

Now applying the robust reformulation introduced above and using the mean and variance estimates from the GPPCE we can reformulate the chance constraints in Equation \ref{eq:path_const_def} and Equation \ref{eq:terminal_const_def} as follows:
\begin{subequations}
\begin{align}
    & \mu_{GP}(\mathbf{z}^{(k)}_{g_j};\bm{\uptau}_{g_j}) + \kappa_{\epsilon} \sqrt{\sigma^2_{GP}(\mathbf{z}^{(k)}_{g_j};\bm{\uptau}_{g_j})} \leq 0 && \forall (k,j) \in \{1,\ldots,N\} \times \{1,\ldots,n_g\}    \\
    & \mu_{GP}(\mathbf{z}_{g_j^N};\bm{\uptau}_{g_j^N}) + \kappa_{\epsilon} \sqrt{\sigma^2_{GP}(\mathbf{z}_{g_j^N};\bm{\uptau}_{g_j^N})} \leq 0 && \forall j \in \{1,\ldots,n_g^N\}
\end{align}
\end{subequations}
where $\mu_{GP}(\mathbf{z}^{(k)}_{g_j};\bm{\uptau}_{g_j})$ and $\mu_{GP}(\mathbf{z}_{g_j^N};\bm{\uptau}_{g_j^N})$ are the GPPCE mean estimates of the constraints $g_j^{(k)}$ and $g_j^N$ respectively based on the data matrices defined in Equation \ref{eq:estimate_data}, while $\sigma^2_{GP}(\mathbf{z}^{(k)}_{g_j};\bm{\uptau}_{g_j})$ and $\sigma^2_{GP}(\mathbf{z}_{g_j^N};\bm{\uptau}_{g_j^N})$ represent the GPPCE variance estimates as introduced in Section \ref{sec:unc_prop_GPPCE}. Note we are using the $GP$ variance that is larger, since it accounts for the error using only \textit{finite number} of samples.      

\subsection{GPPCE SNMPC formulation}
In this section we formulate the stochastic optimal control problem to be solved using the mean and variance approximations of the objective and constraint functions as introduced in Section \ref{sec:uncertainty_prop}. We optimize over the control actions $\mathbf{u}_k$ given the objective and constraints defined in Section \ref{sec:prob_def}:
\begin{equation}
\begin{aligned}
& \underset{\tilde{\mathbf{U}}^{(t)}_{N}}{\text{minimize}} \quad \mu_{GP}(\mathbf{z}_{J^d};\bm{\uptau}_{J^d}) + \omega \cdot \sigma^2_{GP}(\mathbf{z}_{J^d};\bm{\uptau}_{J^d})     \\
& \text{subject to}  \\
& \Tilde{\mathbf{x}}_{k+1}^{(s)} = \mathbf{f}(\Tilde{\mathbf{x}}_k^{(s)},\Tilde{\mathbf{u}}^{(s)}_k,\bm{\uptheta}_s), \quad \Tilde{\mathbf{x}}^{(s)}_k = \hat{\mathbf{x}}_t && \forall (k,s) \in \{t,\ldots,N\} \times \{1,\ldots,n_s\} \\
& \mu_{GP}(\mathbf{z}^{(k)}_{g_j};\bm{\uptau}_{g_j}) + \kappa_{\epsilon} \sqrt{\sigma^2_{GP}(\mathbf{z}^{(k)}_{g_j};\bm{\uptau}_{g_j})} \leq 0 && \forall (k,j) \in \{t,\ldots,N\} \times \{1,\ldots,n_g\} \\
& \mu_{GP}(\mathbf{z}_{g_j^N};\bm{\uptau}_{g_j^N}) + \kappa_{\epsilon} \sqrt{\sigma^2_{GP}(\mathbf{z}_{g_j^N};\bm{\uptau}_{g_j^N})} \leq 0 && \forall j \in \{1,\ldots,n_g^N\} 
\end{aligned}
\label{eq:PCESNMPC}
\end{equation}
where $\hat{\mathbf{x}}_t$ is the current state measurement and $\tilde{\mathbf{U}}^{(t)}_{N}=[\tilde{\mathbf{u}}_t,\ldots,\tilde{\mathbf{u}}_N]^{\sf T}$.

\subsection{GPPCE SNMPC algorithm} \label{sec:GPPCE_alg}
In this section we outline the algorithm to solve the problem defined in Section \ref{sec:prob_def} using the GPPCE SNMPC optimization problem from the previous section. At each time $t$ we are given the current state measurement $\hat{\mathbf{x}}_t$, from which we aim to determine the best control action to take. To formulate the problem the dynamic equation system in Equation \ref{eq:f_x} needs to be defined together with the initial conditions $\hat{\mathbf{x}}_0$ and the time horizon $N$. Further, the objective function $J^d(N,\mathbf{x}_0(\bm{\uptheta}),\bm{\uptheta},\mathbf{U}_N)$ with the variance factor $\omega$ need to be defined, together with the input constraint set $\mathbb{U}$, path constraint functions $g_j(\cdot)$ and terminal constraint functions $g_j^N(\cdot)$. The corresponding probability of feasibility $\epsilon$ needs to be set. Next we specify $\bm{\uptau}_{J^d}$, $\bm{\uptau}_{g_j}$, and $\bm{\uptau}_{g_j^N}$. Lastly, the terms required for the GPPCE estimator are pre-computed according to Equation \ref{eq:pre_computations}. The overall algorithm is stated in Algorithm \ref{alg:GPPCE_SNMPC}.

\begin{algorithm2e}[H] \label{alg:GPPCE_SNMPC}
 \caption{GPPCE SNMPC algorithm}
 \textit{Offline Computations}
 \begin{enumerate}
\item{Choose time horizon $N$, initial condition $\hat{\mathbf{x}}_0$, stage costs $\mathcal{L}(\cdot)$ and terminal cost $\mathcal{M}(\cdot)$, variance weighting factor $\omega$, path constraint functions $g_j(\cdot)$, terminal constraint functions $g_j^N(\cdot)$, input constraint set $\mathbb{U}$, and the chance constraint probability $\epsilon$.}
\item{Specify the GPPCE estimator by setting reasonable values to $\bm{\uptau}_{J^d}$, $\bm{\uptau}_{g_j}$, and $\bm{\uptau}_{g_j^N}$.}
\item{Pre-compute terms required for the GPPCE estimator given in Equation \ref{eq:pre_computations}.}
\end{enumerate} 
\textit{Online Computations} \\ \For{$t=0,\ldots,N-1$}{
\begin{enumerate}
\item{Solve the SNMPC problem in Equation \ref{eq:PCESNMPC} with the current state $\hat{\mathbf{x}}_t$.}
\item{Apply the first control input of the optimal solution to the real plant.} 
\item{Measure the state $\hat{\mathbf{x}}_t$.}
\end{enumerate}}
\end{algorithm2e}

\section{Semi-batch reactor case study} \label{sec:case_study}
\begin{figure}[H]
\centering
\includegraphics[width=0.95\textwidth]{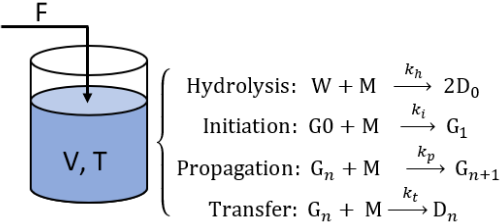}
\caption{Figure summarizing the main variables of the semi-batch reactors with the main reactions taking place. $F$ is the monomer feedrate, $V$ and $T$ are the volume and temperature of the liquid in the reactor respectively, $W$ represents water, $M$ denotes the monomer, $D_n$ and $G_n$ are the dormant and active product chains with length $n$ respectively.} 
\label{fig:semi_batch_reactor}
\end{figure}

The GPPCE SNMPC algorithm introduced in Section \ref{sec:GP_PCE_SNMPC} is applied to a semi-batch reactor case study for the production of the polymer polypropylene glycol from the monomer propylene oxide (PO) as illustrated in Figure \ref{fig:semi_batch_reactor}. 

\subsection{Dynamic model}
For this batch process polymerization reaction an extensive model has been proposed in \citet{Nie2013}. This model has a separate dynamic equation for each chain length separately, which we simplified using the "method of moments" \citep{Nie2013a}. The ordinary differential equations (ODEs) then describe the moments of the polymer as opposed to the amount of each specific chain length. This is often sufficient to estimate important performance criteria. In addition, an energy balance was added due to the importance of temperature control. The dynamic model consists of $7$ ODEs. The dynamic model can be stated as:
\begin{subequations}
\begin{align} \label{eq:case_study_part1}
    & \dot{m}        = F \textit{MW}_{\text{PO}} && m_0 = \hat{m}_0 \\
    & \dot{T}        = \frac{(-\Delta H_p) k_p \gamma_{G0} \textit{M} - \textit{UA}(T-T_C)V - F \textit{MW}_{\text{PO}} C_{pf}(T-T_f)V}{m C_{pb} V}
                       && T_0 = \hat{T}_0 \\ 
    & \dot{W}        = -\frac{k_h W M}{V}    && W_0 = \hat{W}_0           \\
    & \dot{M}        = F - \frac{\left(k_h W + k_i G_0 + k_p \gamma_{G0} + k_t(\gamma_{G0} + G_0)\right) M}{V} && M_0 = \hat{M}_0  \\
    & \dot{X}_0      = \frac{(2 k_h W - k_i G_0)M}{V}      && {X_0}_0 = \hat{X_0}_0        \\
    & \dot{\gamma}_{X0} = \frac{k_i G_0 M}{V}&& {\gamma_{X0}}_0 =  \hat{\gamma_{X0}}_0 \\
    & \dot{\gamma}_{X1} = \frac{(k_i G_0 + k_p \gamma_{G0}) M}{V} && {\gamma_{X1}}_0 = \hat{\gamma_{X1}}_0
\end{align}
\end{subequations}
where $m$ is the liquid mass in $\text{kg}$, $T$ is the reactor temperature in $\text{K}$, $W$ is the amount of water in $\text{kmol}$, $M$ is the amount of monomer in $\text{kmol}$, $X_0$ is the concentration of Methanol in $\text{kmol}$, $\gamma_{X0}$ is the zeroth polymer moment in $\text{kmol}$, and $\gamma_{X1}$ is the first polymer moment in $\text{kmol}$. $F$ is the feed rate of the monomer in $\text{kmol}/\text{s}$ and $T_C$ is the cooling water temperature in $\text{K}$. $\text{k}_\text{p}$, $k_h$, $k_i$, and $k_t$ in $\text{m}^3\text{kmol/s}$ are the kinetic constants of the propagation, hydrolysis, initiation, and transfer reactions respectively. $C_{pb}$ and $C_{pf}$ are the heat capacities of the bulk liquid and the monomer feed respectively in $\text{kJ/kg/K}$. $G_0$, $\gamma_{G0}$, and $\gamma_{G1}$ are the active concentrations of Methanol, the zeroth polymer moment, and the first polymer moment in kmol. $V$ is the liquid volume in the reactor. The kinetic constants and the heat capacities are given as functions of temperature \citep{Nie2013a}:
\begin{subequations}
\begin{align} \label{eq:case_study_part2}
    & k_p    = A_p \exp(-E_{Ap}/RT) \\
    & k_h    = A_h \exp(-E_{Ah}/RT) \\
    & k_i    = A_i \exp(-E_{Ai}/RT) \\
    & k_t    = A_t \exp(-E_{At}/RT) \\
    & C_{pf} = 0.92 + 8.871 \times 10^{-3}  T - 3.1 \times 10^{-5} T^{2} + 4.78 \times 10^{-8} T^3 \\
    & C_{pb} = 1.1 + 2.72 \times 10^{-3} T
\end{align}
\end{subequations}

$G_0$, $\gamma_{G0}$, and $\gamma_{G1}$ depend on $X_0$, $\gamma_{X0}$, and $\gamma_{X1}$ as follows:
\begin{subequations}
\begin{align} \label{eq:case_study_part3}
    & G_0 = X_0 n_C/(X_0 + \gamma_{X0}) \\
    & \gamma_{G0} = \gamma_{X0} n_C/(X_0 + \gamma_{X0}) \\
    & \gamma_{G1} = \gamma_{X1} n_C/(X_0 + \gamma_{X0})
\end{align}
\end{subequations}
where $n_C$ is the amount of catalyst in the reactor in kmol.

Three parameters in the dynamic model were assumed to be uncertain: The pre-exponential factor of the propagation reaction $A_p$ in $\text{m}^3/\text{kmol}/\text{s}$, the overall heat transfer coefficient $\textit{UA}$ in kW/K, and the total amount of catalyst $n_C$ in kmol. The remaining parameters including the initial conditions are given in Table \ref{tab:Parameters}.

\begin{table}[H] \caption{Parameter values for dynamic model taken from \citet{Nie2013} and operating conditions as defined in Equations $\ref{eq:case_study_part1}-\ref{eq:case_study_part3}$.}
\begin{tabular}{llll} \label{tab:Parameters}
Parameter & Value & Units & Description  \\
\hline
$\textit{MW}_{\text{PO}}$ & $58.08$  & $\text{kg/kmol}$ & Molecular weight of PO                   \\
$\Delta H_p$              & $-92048$ & $\text{kJ/kmol}$ & Enthalpy of reaction for propagation reaction                \\
$A_h$                     & $2.4 \times 10^{8}$ & $\text{m}^3/\text{kmol}/\text{s}$ & Pre-exponential factor of hydrolysis kinetic constant \\
$A_i$                     & $4.0 \times 10^{8}$ & $\text{m}^3/\text{kmol}/\text{s}$ & Pre-exponential factor of initiation kinetic constant \\
$A_t$                     & $9.5 \times 10^{8}$ & $\text{m}^3/\text{kmol}/\text{s}$ & Pre-exponential factor of transfer kinetic constant \\
$E_{Ap}$                  & $6.9 \times 10^{4}$  & $\text{kJ/kmol}$ &                  Activation energy of propagation reaction \\
$E_{Ah}$                  & $8.2 \times 10^{4}$  & $\text{kJ/kmol}$ &                  Activation energy of hydrolysis reaction \\
$E_{Ai}$                  & $7.8 \times 10^{4}$ & $\text{kJ/kmol}$ &                  Activation energy of initiation reaction \\
$E_{At}$                  & $1.05 \times 10^{5}$ & $\text{kJ/kmol}$ &                  Activation energy of transfer reaction \\
$R$                       & $8.314$  & $\text{kJ/kmol/K}$ & Universal gas constant \\
$\hat{m}_0$               & $1.56 \times 10^{3}$         & kg   & Initial reactor mass                 \\                    $\hat{T}_0$               & $400$    & K    & Initial reactor temperature          \\          
$\hat{W}_0$               & $1.0$    & kmol & Initial amount of water              \\
$\hat{M}_0$               & $10.0$   & kmol & Initial amount of monomer            \\
$\hat{X_0}_0$             & $0.0$    & kmol & Initial amount of methanol           \\
$\hat{\gamma_{X0}}_0$     & $10.0$   & kmol & Initial zeroth polymer moment        \\
$\hat{\gamma_{X1}}_0$     & $10.0$   & kmol & Initial first polymer moment       
\end{tabular}
\end{table}

The control inputs are given by the monomer feed rate $F$ and cooling water temperature $T_C$. In compact form we can write $\mathbf{x}=[m,T,W,M,X_0,\gamma_{X0},\gamma_{X1}]^{\sf T}$ and $\mathbf{u}=[F,T_C]^{\sf T}$. The uncertain model parameters are given as functions of $\bm{\uptheta}$, which can be used to attain complex probability distributions: 
\begin{subequations}
\begin{align}
    & A_p(\bm{\uptheta}) = \exp\left(-5 + 0.05 \theta_3 + 0.05 \theta_2 + 0.03 \theta_1   \right) \times 10^{9} + 5 \times 10^{6}   \\
    & UA(\bm{\uptheta})  = 40 \cos\left(1.26 + 0.09\theta_3 - 0.09 \theta_2 + 0.09 \theta_1  \right) + 40     \\
    & n_C(\bm{\uptheta}) = |-0.05 \theta_3 + 0.05 \theta_2 - 0.05 \theta_1| + 1
\end{align}
\end{subequations}

\begin{figure}[H]
\centering
\includegraphics[width=1\textwidth]{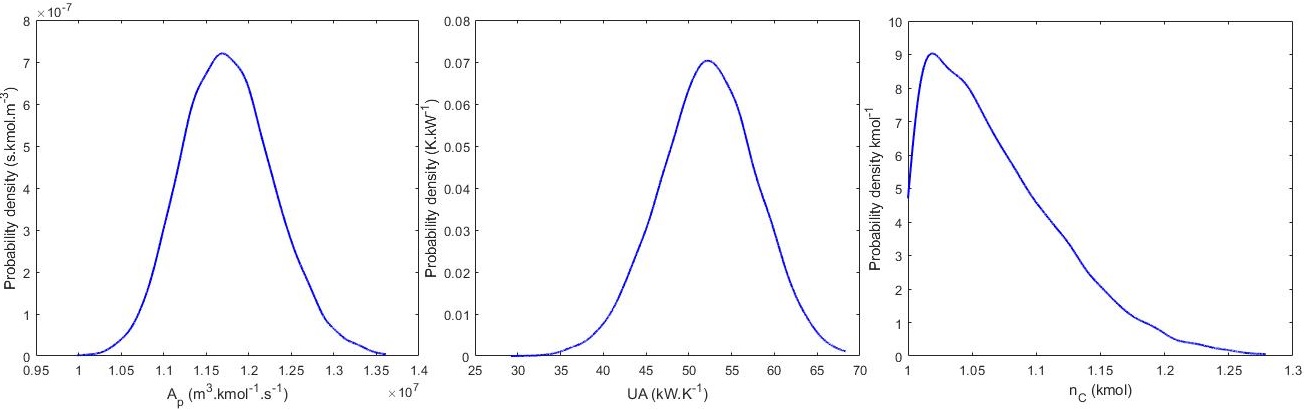}
\caption{Plots of the probability density functions of the uncertain parameters (from left to right for: $A_p$, $\textit{UA}$, $n_C$)} 
\label{fig:Uncertain_parameter_pdfs}
\end{figure}

This allows for the uncertain parameters to attain nearly arbitrary complex probability distributions as can be seen from their respective pdfs plotted in Figure \ref{fig:Uncertain_parameter_pdfs}.

\subsection{Problem set-up}
The time horizon $N$ was set to $12$ with a variable continuous batch time $t_{batch}$ with equal sampling times. The state at each discrete time $t+1$ can be expressed as follows employing Equation \ref{eq:case_study_part1}:
\begin{align}
    & \mathbf{x}_{t+1} = \int_0^{t_{batch}/N} \bar{\mathbf{f}}(\mathbf{x}_t,\mathbf{u}_t) dt + \mathbf{x}_t
\end{align}
where $\bar{\mathbf{f}}(\cdot)=[\dot{m},\dot{T},\dot{W},\dot{M},\dot{X}_0,\dot{\gamma}_{X0},\dot{\gamma}_{X1}]^{\sf T}$ as defined in Equation \ref{eq:case_study_part1}.

The required discrete-time system for Equation \ref{eq:f_x} is then obtained using orthogonal Radau collocation. Each control interval is simulated by a polynomial with an overall degree of $5$. The objective for the control algorithm is aimed to minimize the required batch time $t_{batch}$ with a penalty on changes in the control input:
\begin{align}
    & J^d(N,\hat{\mathbf{x}}_0,\bm{\uptheta},\mathbf{U}_N) = t_{batch} + \sum_{t=1}^N \bm{\Delta}^{\sf T}_{\mathbf{u}_t} \mathbf{R} \bm{\Delta}_{\mathbf{u}_t}
\end{align}
where $\bm{\Delta}_{\mathbf{u}_t} = \mathbf{u}_t - \mathbf{u}_{t-1}$ and $\mathbf{R} = \diag(10^{-6},10^{-4})$. 

The minimization is subject to two terminal constraints and a path constraint. The path constraint aims to keep the reactor temperature below 420K for safety reasons, which can be stated as follows:
\begin{align}
    g(\mathbf{x}_t,\mathbf{u}_t,\bm{\uptheta}) = T - 420 \leq 0
\end{align}

The two terminal constraints state batch product quality properties to be reached. The first terminal constraint requires the batch to reach a number average molecular weight ($\textit{NAMW}$) in kg/kmol of 1500 defined as $\textit{NAMW} = \textit{MW}_{\textit{PO}} \frac{\gamma_{X1}}{\gamma_{X1}}$. The second terminal constraint requires the final monomer concentration to not exceed $1000\text{ppm}$. These terminal constraints are:
\begin{subequations}
\begin{align}
    & g_1^N(\mathbf{x}_N,\bm{\uptheta}) = -\textit{MW}_{\textit{PO}} \frac{\gamma_{X1}}{\gamma_{X0}} + 1500 \leq 0 \\
    & g_2^N(\mathbf{x}_N,\bm{\uptheta}) = 10^6 \times \frac{\textit{MW}_{\textit{PO}} M}{m} - 1000 \leq 0
\end{align}
\end{subequations}
The chance of constraint violation was to $\epsilon = 0.05$ for the constraints defined above. The control inputs are constrained as:
\begin{subequations}
\begin{align}
& 0 \leq F \leq 0.1 \\
& 298.15 \leq T_C \leq 423.15
\end{align}
\end{subequations}

For the GPPCE approximation the scaling variable $\hat{\bm{\uplambda}}$ was set to the following for the different constraint functions, which we refer to as $\textit{anisotropic}$:
\begin{subequations}
\begin{align}
    & \hat{\bm{\uplambda}}_{g}      = [0.22,0.77,0.55] \\
    & \hat{\bm{\uplambda}}_{g_1^N}  = [0.31,0.87,0.44] \\
    & \hat{\bm{\uplambda}}_{g_2^N}  = [0.30,0.75,0.45]
\end{align}
\end{subequations}

The values were determined using Equation \ref{eq:optimal_scaling} by generating different trajectories by setting $\mathbf{U}$ to values between its upper and lower bound. The various $\hat{\bm{\uplambda}}$ values obtained then allowed us to set $\hat{\bm{\uplambda}}$ to reasonable values for the different constraints.    

For comparison purposes we also set the length scales to $\textit{isotropic}$ values using the heuristic in Equation \ref{eq:median_heuristic}, which leads to the following values:  
\begin{subequations} \label{eq:iso_length}
\begin{align}
    & \hat{\bm{\uplambda}}^{iso}_{g} = \hat{\bm{\uplambda}}^{iso}_{g_1^N} = \hat{\bm{\uplambda}}^{iso}_{g_2^N} = [0.65,0.65,0.65]
\end{align}
\end{subequations}

\section{Results and discussions} \label{sec:res_disc}
In this section we present the results of the case study outlined in Section \ref{sec:case_study}. The aim of this section is two-fold. First in Section \ref{sec:GPPCE_accuracy} we compare the accuracy of the GPPCE mean and variance estimates by comparing it to other important approaches that have been utilised to formulate SNMPC problems. Thereafter, in Section \ref{sec:SNMPC_ver} we verify the GPPCE SNMPC algorithm defined in Section \ref{sec:GP_PCE_SNMPC}  and compare it to a nominal SNMPC algorithm with soft constraints. In addition, the GPPPCE SNMPC algorithm is compared to two popular SNMPC algorithms in literature. These are the multi-stage NMPC algorithm \citep{Lucia2014} and a SNMPC formulation using the Unscented transformation \citep{heine2006robust}. 

\subsection{GPPCE accuracy} \label{sec:GPPCE_accuracy}
\begin{figure}[H]
\centering
\includegraphics[width=1\textwidth]{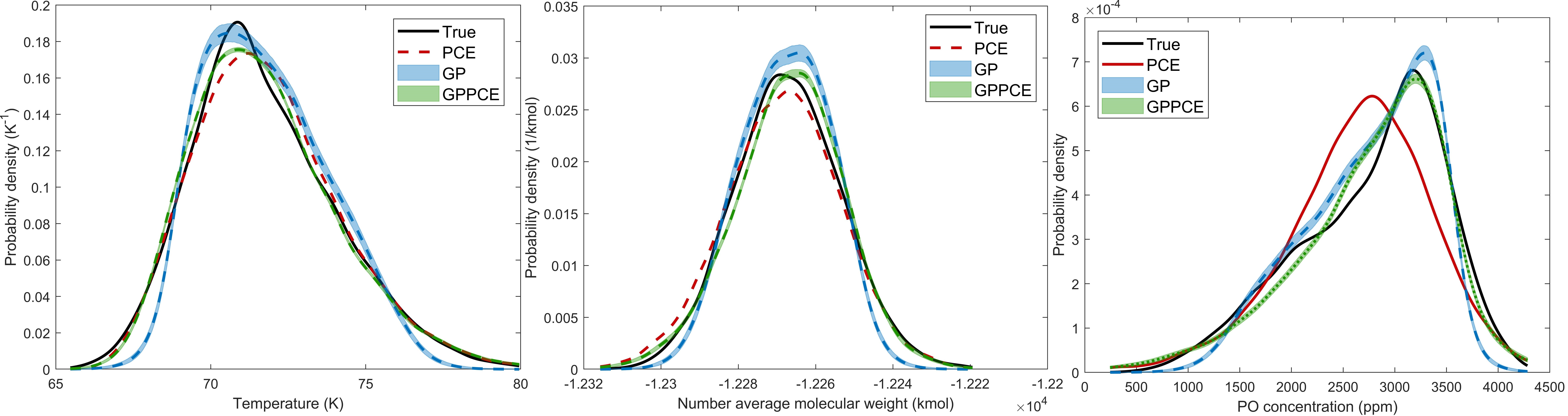}
\caption{Plots of the probability density functions of the models $PCE$, $GP$, and $GPPCE$ for the path constraint function $g(\cdot)$ at $t=N$ and the two terminal constraint functions $g_1^N(\cdot)$, $g_2^N(\cdot)$ from left to right respectively for $\mathbf{U}$ set to its upper bound.} 
\label{fig:GP_PCE_GPPCE_pdfs}
\end{figure}

\begin{figure}[H]
\centering
\includegraphics[width=1\textwidth]{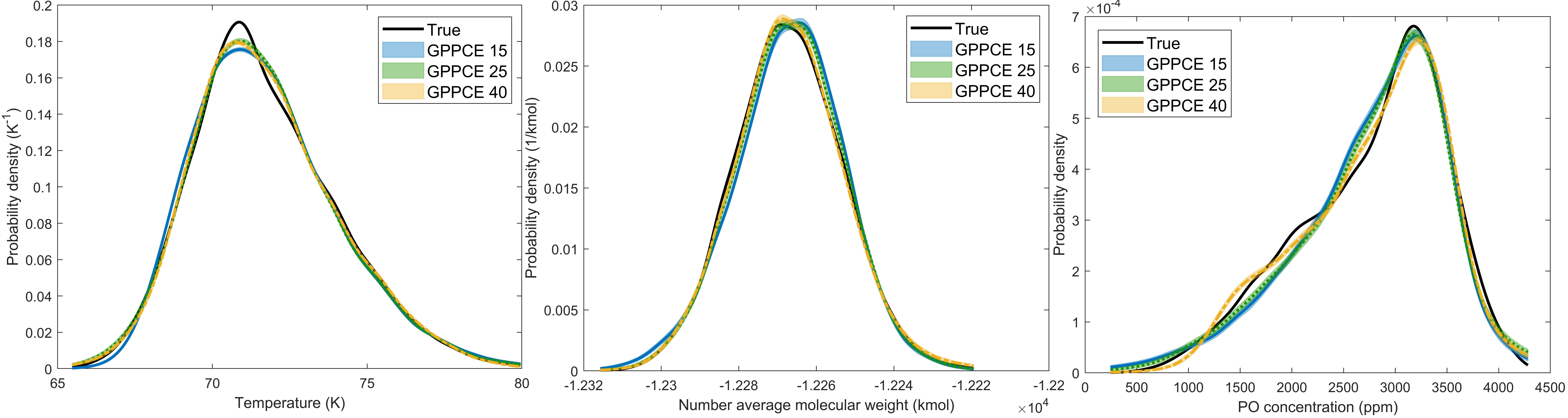}
\caption{Plots of the probability density functions of the $GPPCE$ model with $15$, $25$, and $40$ training data-points for the path constraint function $g(\cdot)$ at $t=N$ and the two terminal constraint functions $g_1^N(\cdot)$, $g_2^N(\cdot)$ from left to right respectively for $\mathbf{U}$ set to its upper bound.} 
\label{fig:GPPCE15_GPPCE25_GPPCE40_pdfs}
\end{figure}

\begin{figure}[H]
\centering
\includegraphics[width=0.85\textwidth]{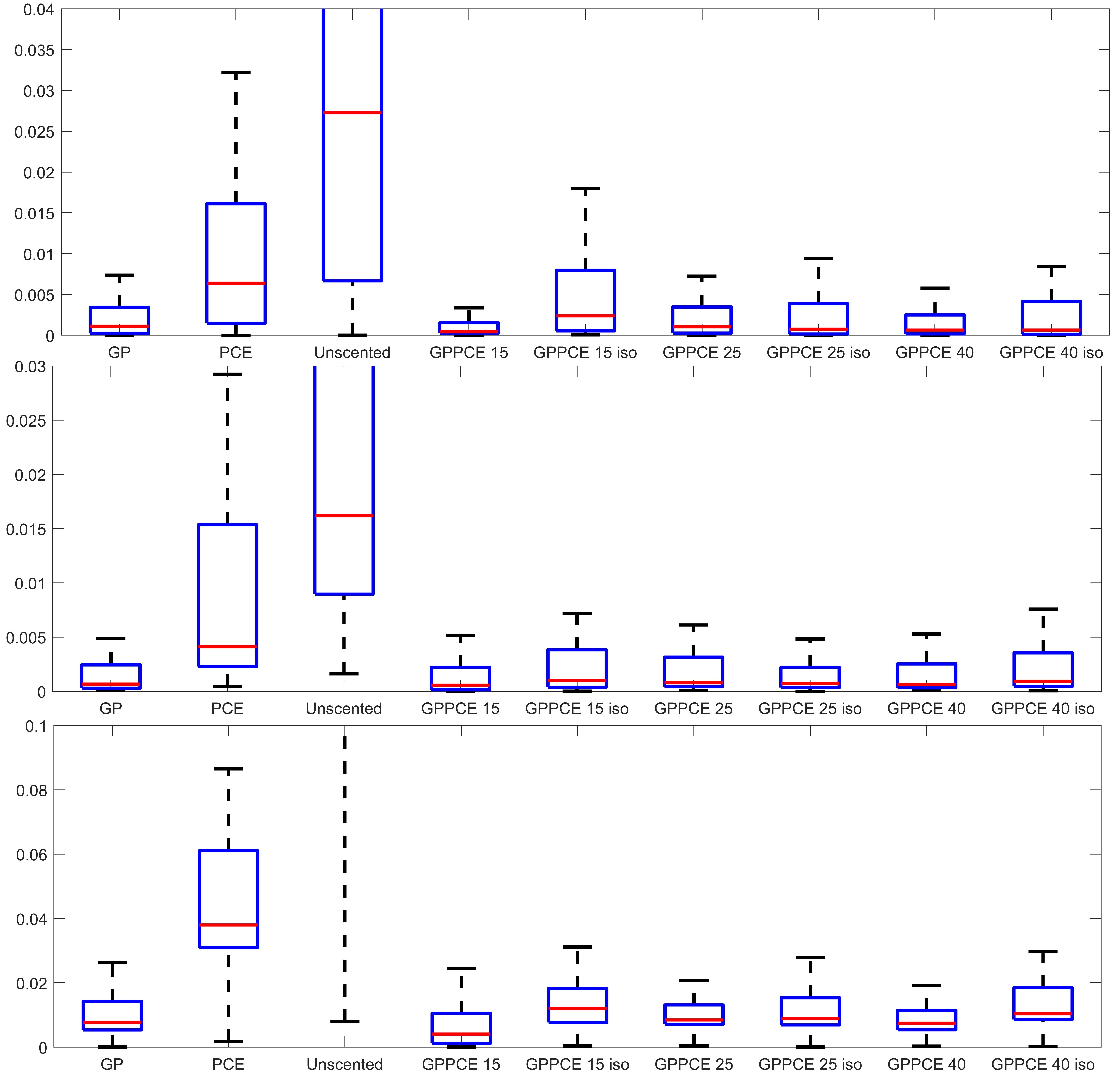}
\caption{Box plots of absolute relative error of mean estimates from 100 random $\mathbf{U}$ values. From top to bottom the values in each field refer to the path constraint $g(\cdot)$ at $t=N$, the first terminal constraint $g_1^N(\cdot)$, and the second terminal constraint $g_2^N(\cdot)$ respectively.} 
\label{fig:Box_mean_error}
\end{figure}

\begin{figure}[H]
\centering
\includegraphics[width=0.85\textwidth]{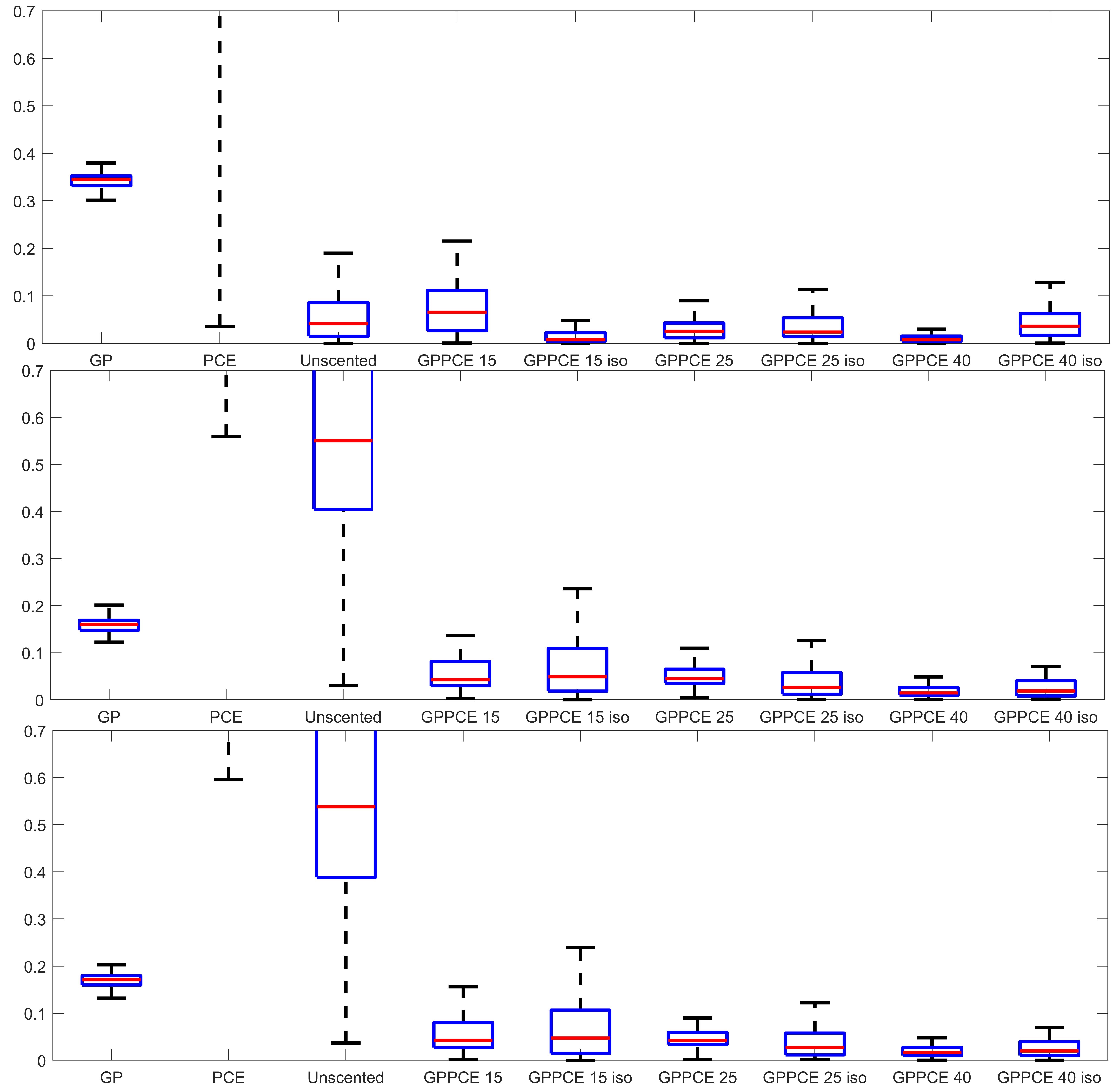}
\caption{Box plots of absolute relative error of standard deviation estimates from 100 random $\mathbf{U}$ values. From top to bottom the values in each field refer to the path constraint $g(\cdot)$ at $t=N$, the first terminal constraint $g_1^N(\cdot)$, and the second terminal constraint $g_2^N(\cdot)$ respectively.} 
\label{fig:Box_var_error}
\end{figure}

In this section we verify the GPPCE approach to obtain mean and variance estimates of nonlinear transformations as outlined in Section \ref{sec:GP_PCE}. To accomplish this we ran the following tests:
\begin{itemize}
    \item{Set $\mathbf{U}$ to its upper bound and compare the pdfs of PCE, GP, and GPPCE for the two terminal constraint functions and the path constraint function at the final time $t=N$ with the \textit{true} pdfs. Each model has $15$ training data points according to a Sobol design and the PCE and GPPCE polynomial order was set to $2$. Note the pdfs are obtained using kernel density estimation (KDE) of the respective models \citep{Silverman2018}. The models are obtained using the data and polynomial order as outlined in the previous section. The results for this are shown in Figure \ref{fig:GP_PCE_GPPCE_pdfs} with the corresponding mean and variance estimates given in Table \ref{tab:Variances_means}.}
    \item{Set $\mathbf{U}$ to its upper bound and compare the pdfs of GPPCE with $15$, $25$, and $40$ training data-points according to a Sobol design with polynomial order of $2$ throughout for the two terminal constraint functions and the path constraint function at the final time $t=N$. The results for this are shown in Figure \ref{fig:GPPCE15_GPPCE25_GPPCE40_pdfs}. Furthermore, the same graphs are shown for the case with isotropic (iso) length scales given in Equation \ref{eq:iso_length} for GPPCE with $15$, $25$, and $40$ training data-points.}
    \item{Lastly, $\mathbf{U}$ was set to $100$ random values. For these the mean and variance of the terminal constraint functions and the path constraint at $t=N$ were estimated using GPs with $15$ data-points, PCEs with $15$ data-points, the Unscented transformation with $7$ data-points, and GPPCE based on $15$, $25$, and $40$ data-points using anisotropic length scales and isotropic length scales. The Unscented transformation has a fixed number of data-points, which corresponds to two times the number of inputs plus one (2$n_{\bm{\uptheta}}$ + 1). The relative absolute error for these is illustrated in Figure $\ref{fig:Box_mean_error}$ for the mean estimates  and in Figure $\ref{fig:Box_var_error}$ for the standard deviation estimates as box plots.}
\end{itemize}

\begin{figure}[H]
\centering
\includegraphics[width=1\textwidth]{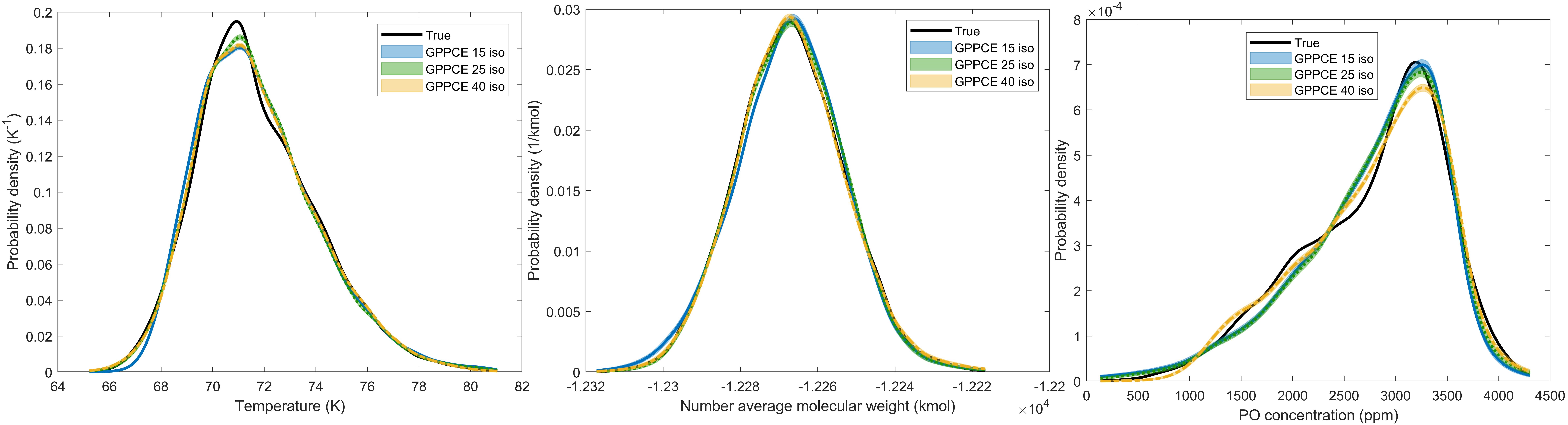}
\caption{Plots of the probability density functions of the $GPPCE$ model with $15$, $25$, and $40$ training data-points with isotropic length scales for the path constraint function $g(\cdot)$ at $t=N$ and the two terminal constraint functions $g_1^N(\cdot)$, $g_2^N(\cdot)$ from left to right respectively for $\mathbf{U}$ set to its upper bound.} 
\label{fig:GPPCE15_GPPCE25_GPPCE40_pdfs_iso}
\end{figure}

Based on the tests outlined above we can draw the following observations and conclusions:
\begin{itemize}
    \item{Generally speaking from the plots in Figure \ref{fig:GP_PCE_GPPCE_pdfs} and Table \ref{tab:Variances_means} it can be said that all three approaches are able to represent the pdfs reasonably well with good approximations to the mean and variance. The mean values are accurate within a $4\%$ error margin for all approaches, where the PCE has the largest error for $g_2^N(\cdot)$. The variance is generally more difficult to approximate with an error margin of $40\%$ due to inaccuracies of the GP approximation. PCEs are seen to outperform GPs considerably for the mid-plot, while for the graph on the RHS GPs outperform PCEs. GPs are expected to be able to handle more complex responses due to interpolating between data-points, while PCEs often capture better the overall trend, i.e. lead to a better global fit. GPPCE seems to approximate both well, since it is based on both methods, which highlights its major advantage over both. For example GPPCE 15 has a percentage variance error of at most $10\%$, compared to PCE with $15\%$ and GP with $40\%$. Lastly, comparing Figure \ref{fig:GPPCE15_GPPCE25_GPPCE40_pdfs} and Figure \ref{fig:GPPCE15_GPPCE25_GPPCE40_pdfs_iso} it can be seen that the anisotropic length-scales leads to a better fit for the pdfs. From Table \ref{tab:Variances_means} the maximum error of the variance for the isotropic case is around $13\%$ compared to only $9.5\%$ in the anisotropic case.}
    \item{The confidence bound for GPPCEs and GPs in Figure \ref{fig:GP_PCE_GPPCE_pdfs} and Figure \ref{fig:GPPCE15_GPPCE25_GPPCE40_pdfs} corresponds to a $95\%$ confidence region. It can be seen that the region is able to capture the relative uncertainty well, since it is larger if the fit is poor and smaller if the fit is better. Nonetheless, it does seem to be often overconfident. Furthermore, the variance accounting for the \textit{finite samples} is given in Table \ref{tab:Variances_means} as \textit{stochastic variance}. It can be seen that this variance is only once smaller then the \textit{true} variance, further highlighting its potential use as a more conservative variance estimate. In addition, the tails of the pdfs are not well captured by the confidence region of the GP. This is expected, since the GP uncertainty measure in Equation \ref{eq:GP_posterior} is upper bounded by $\hat{\alpha}$ and simply assumes its maximum in these regions due to sparsity of data.}
    \item{In Figure \ref{fig:Box_mean_error} the box plots highlight the absolute error from $100$ mean approximations. It can be seen that in general the Unscented transformation performs the worst, which is not surprising since it is based on \textit{only} 7 data-points. Furthermore, PCEs are seen to perform rather poorly as well. For PCEs it should be noted that for several variations it performed very well as can be seen in Figure \ref{fig:GP_PCE_GPPCE_pdfs} for example, however for several variations it performed poorly. This can in part be explained by the difficulty of determining reasonable regularization parameters, which can be viewed as a significant disadvantage. GPs on the other hand often perform worse than PCEs, but manage to never perform very poorly due their nature of interpolating between data-points. Therefore, GPs on average perform much better than PCEs as can be seen in Figure \ref{fig:Box_mean_error}. Lastly, GPPCE can again be seen to outperform both PCEs and GPs, which is in line with previous observations that GPPCE captures the \textit{best} of both techniques. Interestingly the GPPCE mean approximation does not seem to improve with more data-points, however its \textit{worst} performance is already small at $3\%$ with 15 data-points.}
    \item{In Figure \ref{fig:Box_var_error} the box plots show the absolute error from $100$ variance approximations. It should be noted that accurate variance estimates are more difficult to achieve. The PCE performs the worst in this case, because as previously it is performing very poorly on a few test cases. This is further exacerbated from the square variance definition. Unscented performs poorly again, except for the second terminal constraint. GPs also do not perform particularly well leading to a up to nearly $40\%$ error for the second terminal constraint. GPPCE on the other hand performs much better with 15 data-points leading to an error of at most $15\%$ for the first terminal constraint. Further, it can be seen that GPPCE variance approximations steadily improves with more data-points, with 40 data-point GPPCE never exceeding a $5\%$ error threshold. Overall it can be said that GPPCE is a vast improvement over its GP and PCE counterparts for estimating variances.}
    \item{Lastly, the GPPCE with isotropic length-scales performs nearly always worse than using anisotropic length-scales in 8 out of 9 cases for the means, and 8 out of 9 cases for the variances. It should be noted however that isotropic length scales still perform well and are still superior to the GP, PCE, and Unscented approximations for all variance estimates.}
\end{itemize}

\begin{table}[H] \caption{Mean, variance, and \textit{stochastic} variance for the PCE, GP, GPPCE with 15, 25 and 40 training data-points. From left to right the values in each field refer to the path constraint $g(\cdot)$ at $t=N$, the first terminal constraint $g_1^N(\cdot)$, and the second terminal constraint $g_2^N(\cdot)$ respectively.}

\begin{tabular}{llll} \label{tab:Variances_means}
Estimator & Mean  & Variance  & Stochastic variance  \\
\hline
True                      & $71.7, \; -1.23 \times 10^{4}, \; 2780$  & $5.50, \; 195.1, \; 4.96 \times 10^{5}$ &                    \\
PCE                       & $71.9, \; -1.23 \times 10^{4}, \; 2660$  & $5.59, \; 222.0, \; 4.96 \times 10^{5}$ &                    \\
GP                        & $71.8, \; -1.23 \times 10^{4}, \; 2760$  & $3.57, \; 118.5, \; 3.40 \times 10^{5}$ & $4.01, \; 203.4, \; 5.73 \times 10^{5}$                    \\
GPPCE 15                  & $71.8, \; -1.23 \times 10^{4}, \; 2778$  & $5.60, \; 199.1, \; 5.43 \times 10^{5}$ & $5.67, \; 203.4, \; 5.73 \times 10^{5}$ \\
GPPCE 25                  & $71.7, \; -1.23 \times 10^{4}, \; 2809$  & $5.67, \; 189.2, \; 5.51 \times 10^{5}$ & $5.74, \; 196.8, \; 5.85 \times 10^{5}$ \\
GPPCE 40                  & $71.7, \; -1.23 \times 10^{4}, \; 2803$  & $5.50, \; 202.0, \; 5.12 \times 10^{5}$ & $5.55, \; 204.4, \; 5.33 \times 10^{5}$ \\
GPPCE 15 iso                  & $71.8, \; -1.23 \times 10^{4}, \; 2759$  & $5.66, \; 202.3, \; 5.56 \times 10^{5}$ & $5.77, \; 207.2, \; 5.82 \times 10^{5}$ \\
GPPCE 25 iso                  & $71.7, \; -1.23 \times 10^{4}, \; 2795$  & $5.64, \; 188.4, \; 5.61 \times 10^{5}$ & $5.83, \; 192.8, \; 5.96 \times 10^{5}$ \\
GPPCE 40 iso                  & $71.7, \; -1.23 \times 10^{4}, \; 2798$  & $5.60, \; 203.5, \; 5.22 \times 10^{5}$ & $5.71, \; 205.2, \; 5.35 \times 10^{5}$
\end{tabular}
\end{table}

\subsection{SNMPC verification} \label{sec:SNMPC_ver}
To verify the SNMPC algorithm given in Section \ref{sec:GPPCE_alg} we run 400 closed-loop MC simulations of the case study outlined in Section \ref{sec:case_study} by sampling the uncertain parameters $\bm{\uptheta}$ independently. For comparison purposes the $400$ MC simulations of the SNMPC algorithm are compared to $400$ MC simulations of a nominal NMPC algorithm with soft constraints. In addition, the GPPPCE SNMPC algorithm is compared to the the multi-stage NMPC algorithm \citep{Lucia2014} and a SNMPC formulation using the Unscented transformation \citep{heine2006robust}. These algorithms are as follows:
\begin{itemize}
\item{SNMPC: The GPPCE SNMPC algorithm outlined in Section \ref{sec:GPPCE_alg} with 15 data-points.}
\item{Nominal NMPC: A nominal NMPC algorithm with soft constraints ignoring the parametric uncertainties.}
\item{Multi-stage: The multi-stage NMPC algorithm proposed in \citet{Lucia2014} that employs independent scenarios of the parametric uncertainties with separate control actions. We used 15 realizations of the parametric uncertainties that were chosen to give a relatively large spread. }
\item{Unscented: A SNMPC algorithm using the Unscented transformation to obtain the mean and variance estimates to evaluate the probabilistic objective and chance constraints. This should highlight the importance of accurate mean and variance estimates given the poor performance of the Unscented transformation in Section \ref{sec:GPPCE_accuracy}.}
\end{itemize}

\begin{figure}[H]
\centering
\includegraphics[width=0.95\textwidth]{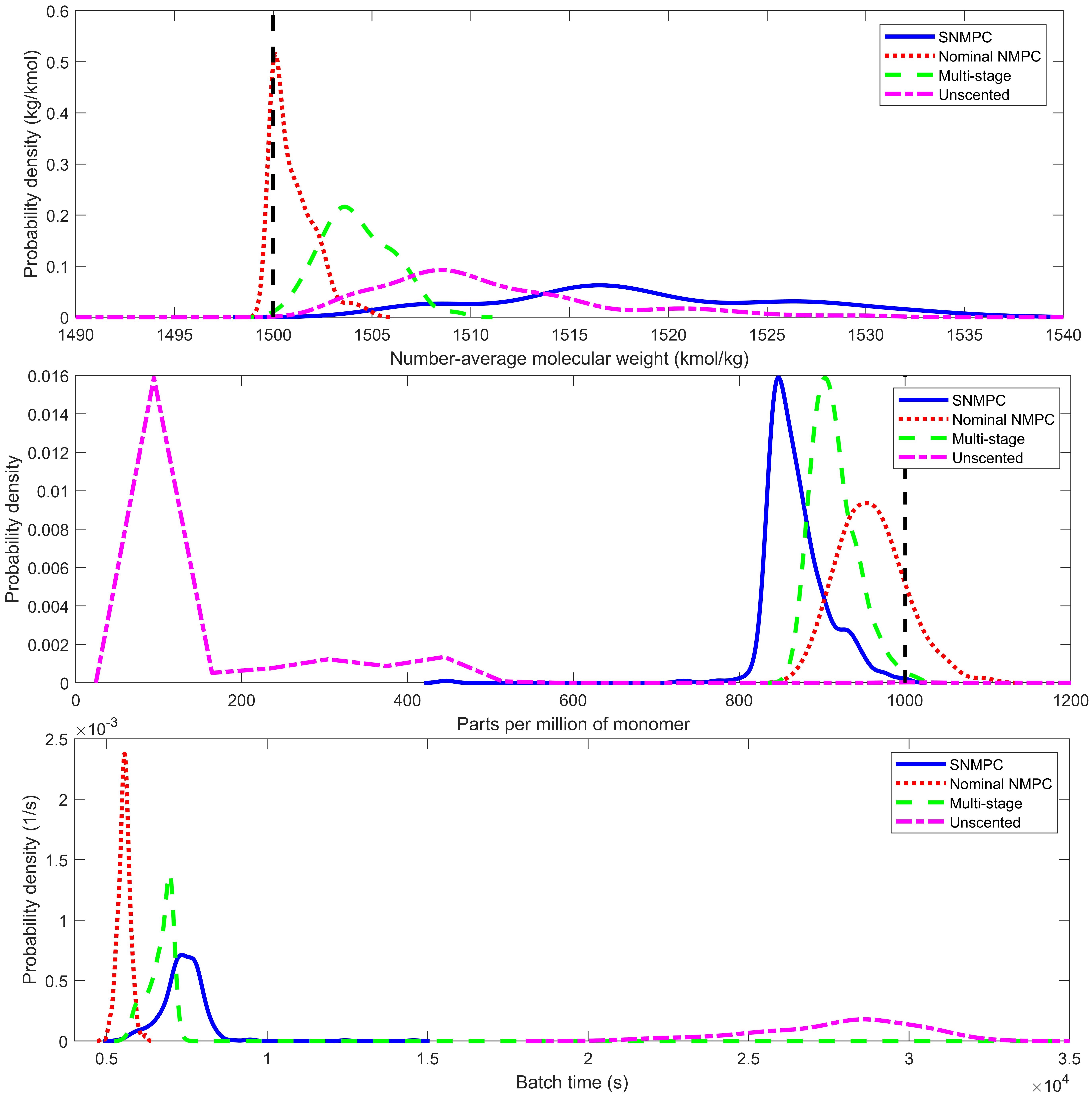}
\caption{Plots of probability densities of \textit{NAMW} (kg/kmol), parts per million of monomer, and batch time (s) based on $400$ MC simulations from top to bottom respectively.} 
\label{fig:MC_pdfs}
\end{figure}

\begin{figure}[H]
\centering
\includegraphics[width=0.95\textwidth]{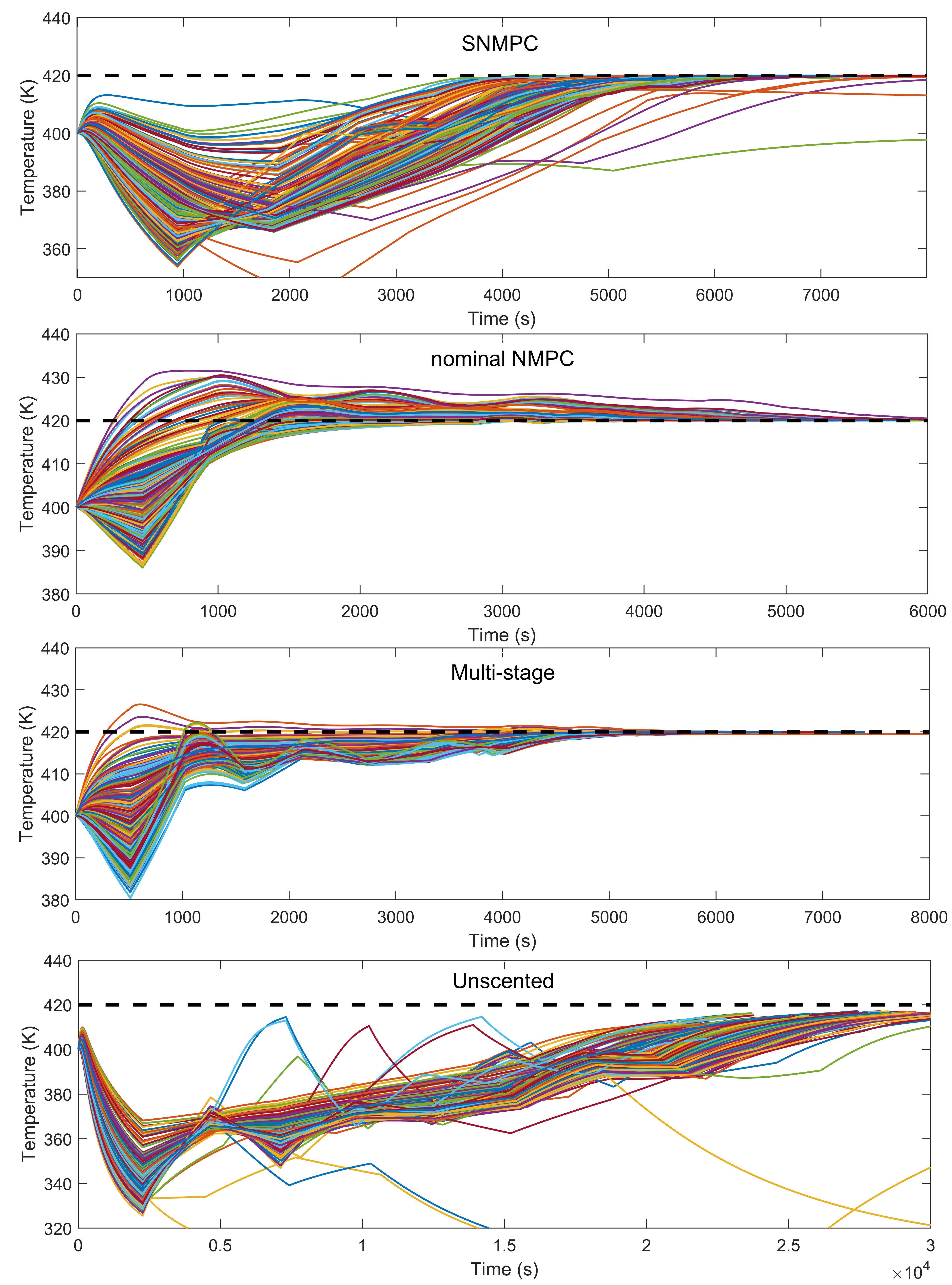}
\caption{Plots of temperature trajectories of $400$ MC simulations.} 
\label{fig:SNMPC_temperatures}
\end{figure}

The results of these MC simulations are highlighted in Figure \ref{fig:MC_pdfs} and Figure \ref{fig:SNMPC_temperatures}. Figure \ref{fig:MC_pdfs} depicts the probability densities using KDE of the 400 MC simulations for the two terminal constraint functions and the final batch time required. Figure \ref{fig:SNMPC_temperatures} shows the temperature trajectories of all $400$ MC simulations for each algorithm. In Table \ref{tab:OCP_time} the average and standard deviation computational times are given. Based on these results we can draw the following observations and conclusions:

\begin{itemize}
    \item{In Figure \ref{fig:MC_pdfs} on the top we can see from the pdf the nominal NMPC violates the constraint considerably more than the SNMPC algorithm reaching often not the required \textit{NAMW}. For all $400$ MC simulations the SNMPC algorithm does not violate this constraint even once, while the nominal NMPC algorithm violates the constraint in $70\%$ of the simulations. The multi-stage NMPC algorithm violates this constraint only once in the $400$ runs, while the Unscented NMPC never violates it. The multi-stage NMPC can be seen to be less conservative staying closer to the constraint on  average.}
    \item{In Figure \ref{fig:MC_pdfs} the mid plot shows the pdfs of the ppm of the monomer. Again it can be seen that the nominal NMPC violates this constraint frequently, while the SNMPC is considerably more robust. Again the SNMPC algorithm does not violate this constraint for all $400$ MC simulations, while the nominal NMPC algorithm violates it $62\%$ of the time. The multi-stage NMPC violates the constraint 3 times in the $400$ runs, while the Unscented NMPC does not violate it at all. The multi-stage NMPC does manage to stay closer to the constraint and can therefore be said to be less conservative, while for the Unscented NMPC the constraint seems not to be active for nearly all scenarios.}
    \item{From Figure \ref{fig:SNMPC_temperatures} it can be seen that the temperature control of the SNMPC algorithm is considerably improved over the nominal NMPC, which violates the constraint in many of the MC simulations. The SNMPC was found to not violate the temperature constraint at all, while the nominal NMPC algorithm violates the path constraint in $95\%$ of cases. The multi-stage NMPC violates the temperature constraint quite substantially in $16\%$ of the scenarios with a maximum violation of up to $8$K. The Unscented NMPC on the other hand manages to not violate the temperature constraint for all scenarios, however the trajectories are very conservative. This may be the result of a severe overestimation of the true variance of the temperature constraint, which could be expected given the previous results in Section \ref{sec:GPPCE_accuracy}.}
    \item{The bottom plot in Figure \ref{fig:MC_pdfs} shows the trade-off of the improved constraint satisfaction for the SNMPC algorithm. The pdf of the batch time for the SNMPC algorithm is more skewed towards longer batch times to be able to adhere the constraints. While the nominal NMPC algorithm on average has a batch time of $5550$s, the SNMPC algorithm takes on average $7380$s to complete a batch. This is expected, since superior constraint satisfaction leads to a worse objective. The multi-stage NMPC algorithm is less conservative than the SNMPC algorithm requiring on average $6900$s, however this may be the result of violating the temperature constraints. The Unscented NMPC has significantly higher batch times of $28,000$s, which are nearly $4$ times as high as the SNMPC approach.}
    \item{Table \ref{tab:OCP_time} shows the computational times of the algorithms. It can be seen that the SNMPC algorithm has an average computational times of $2.9$s, which is $63$ times longer than the computational times of the nominal NMPC algorithm. Nominal NMPC is based on a much smaller optimization problem without scenarios and has less strict constraints, which are often easier to adhere. Also the absence of hard constraints can lead to considerably faster computational times. The standard deviation is very large compared to the mean value for the SNMPC, which indicates a highly asymmetric distribution of computational times. This is found to be due to some OCP times to be considerably larger in excess of 20 seconds, which happens in approximately $2\%$ of the scenarios. The multi-stage NMPC has more than $50\%$ higher average computational times compared to the SNMPC method at $5.4$s due to the increased number of decision variables, however the standard deviation is considerably smaller. The Unscented algorithm has the highest average computational times with $5.4$s, which is nearly double the computational times of the SNMPC approach. The standard deviation is also the highest at $12.2$s. This is most likely due to the large batch times and high conservativenss of the optimization problem.}
    \item{In the presented example none of the constraints are violated for the SNMPC algorithm over 400 closed loop runs, while the percentage constraint violation was set to 0.05. This is very conservative and a result of Chebyshev's inequality, since for a probability of $0.05$ we might expect some constraints to be violated for around $20$ of the 400 scenarios.}
    \item{All in all, it has been shown that the SNMPC algorithm is able to adhere the constraints despite the stochastic parametric uncertainties, while the multi-stage NMPC violated the temperature constraints substantially. This may be due to the multi-stage NMPC not being able to account for intermediary uncertainty values, since it is based entirely on discrete realizations of the uncertainty. Further, it can be seen from the Unscented transformation based NMPC that accurate mean and variance estimates are vital. In Section \ref{sec:GPPCE_accuracy} it was shown that GPPCE is able to attain accurate mean and variance estimates for this case study, while the Unscented transformation leads to large errors. The Unscented NMPC could be seen to consistently overestimate the variances of the constraints, which leads to a very conservative solution and consequently to batch times that are on average $4$ times longer than for the SNMPC approach. Lastly, it is paramount to account for the uncertainty, since otherwise large constraint violations are inevitable as seen by the nominal NMPC method.}
\end{itemize}

\begin{table}[H]
\centering
\caption{Optimal control problem average and standard deviation of computational times.}
\begin{tabular}{*{3}{l}}
Algorithm & average OCP time (s) & standard deviation OCP time (s)  \\
\hline
SNMPC        & 2.868 & 9.277 \\
Nominal NMPC & 0.045 & 0.001 \\
Unscented    & 5.449 & 12.20 \\
Multi-stage  & 4.562 & 0.154
\end{tabular}
\label{tab:OCP_time}
\end{table}

\section{Conclusions} \label{sec:conclusions}
In conclusion, we proposed a new approach to approximate the mean and variance of a nonlinear transformation given a standard normally distributed input by combining GPs and PCEs. It was shown that the method in general is able to capture better the shape of pdfs and leads to much improved approximations of both mean and variance. This can in part be explained by the approach leading to a good approximation if either GP or PCE lead to a good fit. Further, the GPPCE SNMPC algorithm is shown to lead to superior constraint satisfaction over a nominal NMPC algorithm using soft constraints and the multi-stage NMPC approach despite the stochastic uncertainties present. Further, an Unscented transformation based NMPC algorithm required on average $4$ times longer batch times due to poor variance estimates. This highlights the importance of accurate mean and variance estimates for the SNMPC algorithm to work well. The computational times are kept moderately low by pre-computing the expensive terms involved in the GPPCE approach.    

\section*{Acknowledgements}
This project has received funding from the European Union's Horizon 2020 research and innovation programme under the Marie Sklodowska-Curie grant agreement No 675215.

\bibliography{GP_based_MPC}

\appendix \label{sec:expectation_derivation}
\section{Posterior mean and variance derivation}
In Section \ref{sec:GP_PCE} we show how to derive the posterior mean and variance given a fitted GPPCE. The covariance function is given by the SE covariance function as given in Equation \ref{eq:SE} and for the mean function we use Hermite polynomials, as is done in PCE \citep{Kersaudy2015}. Given these choices we need to determine several expectations for Equation \ref{eq:pre_computations} with $\bm{\uptheta} \sim \mathcal{N}(\bm{\uptheta};\mathbf{0},\mathbf{I})$ following a standard normal distribution. Given these choices we can derive the expectations we require in turn. In essence we will exploit the fact that the SE covariance function represents an unnormalized Gaussian pdf:
\begin{equation} \label{eq:SE_normal}
    k(\bm{\uptheta},\bm{\uptheta}_i) = \frac{\hat{\alpha}^2 \mathcal{N}(\bm{\uptheta};\bm{\uptheta}_i,\hat{\bm{\Lambda}})}{|2\pi \hat{\bm{\Lambda}}|^{-\frac{1}{2}}}
\end{equation}
and exploit the following identity:
\begin{align}
    \mathcal{N}(\bm{\uptheta};\bm{\upmu}_{p},\bm{\Sigma}_{p}) \cdot \mathcal{N}(\bm{\uptheta};\bm{\upmu}_{q},\bm{\Sigma}_{q}) \coloneqq r \mathcal{N}(\bm{\uptheta};\bm{\upmu}_{r},\bm{\Sigma}_{r}) \\
    \int \mathcal{N}(\bm{\uptheta};\mathbf{0},\mathbf{I}) \exp\left(-\frac{1}{2} \bm{\uptheta}^{\sf T} \bm{\Lambda}^{-1} \bm{\uptheta}   \right) d\bm{\uptheta} \coloneqq |\mathbf{I} + \bm{\Lambda}^{-1}|^{-\frac{1}{2}}
\end{align}
where $\bm{\Sigma}_{r} = (\bm{\Sigma}_{p}^{-1} + \bm{\Sigma}_{q}^{-1})^{-1}$, $\bm{\upmu}_{r} = \bm{\Sigma}_{r}\left(\bm{\Sigma}_{p}^{-1}\bm{\upmu}_{p} + \bm{\Sigma}_{q}^{-1}\bm{\upmu}_{q} \right)$, and $r = |2 \pi (\bm{\Sigma}_{p} + \bm{\Sigma}_{q})|^{-\frac{1}{2}} \exp \left(-\frac{1}{2} (\bm{\upmu}_{p}-\bm{\upmu}_{q})^{\sf T} (\bm{\Sigma}_{p} + \bm{\Sigma}_{q})^{-1}(\bm{\upmu}_{p}-\bm{\upmu}_{q})\right)$. 

\subsection*{Expectation of $m(\bm{\uptheta})$}
The expectation of $m(\bm{\uptheta})$ is given by the first expansion coefficient in Equation \ref{eq:mean_fcn} due to the orthogonality properties of the Hermite polynomials \citep{Mesbah2014}:
\begin{equation}
\mu_m = \hat{\beta}_0   
\end{equation}
where $\mu_m = \mathbb{E}_{\bm{\uptheta}}[m(\bm{\uptheta})]$.

\subsection*{Expectation of  $ \mathbf{k}_{\zeta,\mathbf{z}}(\bm{\uptheta})$}
The expectation of $\mathbf{k}_{\zeta,\mathbf{z}}(\bm{\uptheta})$ can be derived by expressing the SE covariance function as a multivariate normal distribution. The full derivation can be found in \citet{Deisenroth2011} and leads to the following: 
\begin{align}
    & [{\mu_{\mathbf{k}_{\zeta,\mathbf{z}}}}]_{i} = \hat{\alpha}^2 |\mathbf{I}+\bm{\Lambda}^{-1}|^{-\frac{1}{2}} \exp \left(-\frac{1}{2} \bm{\uptheta}_i^{\sf T}  (\mathbf{I} + \bm{\Lambda})^{-1} \bm{\uptheta}_i \right) 
\end{align}
where ${\bm{\upmu}_{\mathbf{k}_{\zeta,\mathbf{z}}}} = \mathbb{E}_{\bm{\uptheta}}\left[\mathbf{k}_{\zeta,\mathbf{z}}(\bm{\uptheta})\right]$.

\subsection*{Expectation of  $\left(m(\bm{\uptheta})\right)^2$}
Note this is by definition the second moment of $m(\bm{\uptheta})$ and hence has received considerable attention. Again due to the orthogonality properties of the Hermite polynomials utilized the expectation of $\left(m(\bm{\uptheta})\right)^2$ is considerably simplified \citep{Mesbah2014}:
\begin{align}
\mu_{m^2} = \sum_{i=0}^{L-1} \beta_i^2 \mathbb{E}_{\bm{\uptheta}}\left[\phi_i^2(\bm{\uptheta}) \right]
\end{align}
where $\mu_{m^2} = \mathbb{E}_{\bm{\uptheta}}\left[m(\bm{\uptheta})^2\right]$. 

\subsection*{Expectation of $\mathbf{k}_{\zeta,\mathbf{z}}(\bm{\uptheta}) \mathbf{k}_{\zeta,\mathbf{z}}^{\sf T}(\bm{\uptheta})$}
The expectation of the outer product of $\mathbf{k}_{\zeta,\mathbf{z}}(\bm{\uptheta})$ can be found in \citet{Deisenroth2011} and is as follows: 
\begin{align}
    [\mathbf{M}_{\mathbf{k}_{\zeta,\mathbf{z}} \mathbf{k}_{\zeta,\mathbf{z}}^{\sf T}}]_{ij} = k(\bm{\uptheta}_i,\mathbf{0})k(\bm{\uptheta}_j,\mathbf{0})|\mathbf{R}|^{-\frac{1}{2}} \exp\left(\mathbf{l}^{\sf T} \mathbf{R}^{-1} \mathbf{l} \right)
\end{align}
where $\mathbf{M}_{\mathbf{k}_{\zeta,\mathbf{z}} \mathbf{k}_{\zeta,\mathbf{z}}^{\sf T}} = \mathbb{E}\left[    \mathbf{k}_{\zeta,\mathbf{z}}(\bm{\uptheta}) \mathbf{k}_{\zeta,\mathbf{z}}^{\sf T}(\bm{\uptheta}) \right]$, $\mathbf{R}=2\hat{\bm{\Lambda}}^{-1} + \mathbf{I}$, and $\mathbf{l}=\hat{\bm{\Lambda}} \bm{\uptheta}_i + \hat{\bm{\Lambda}} \bm{\uptheta}_j$.

\subsection*{Expectation of $\mathbb{E}_{\bm{\uptheta}}\left[m(\theta) \mathbf{k}_{\zeta,\mathbf{z}}(\bm{\uptheta})\right]$}
This term is somewhat more difficult to deal with, since it is a cross-term between the mean function and the covariance function. Unfortunately we cannot exploit the orthogonality properties of $m(\bm{\uptheta})$, nonetheless the term has a closed-form solution as we will show here:
\begin{align}
    & \left[\bm{\upmu}_{m \mathbf{k}_{\zeta,\mathbf{z}}} \right]_i = \int \mathcal{N}(\bm{\uptheta};\mathbf{0},\mathbf{I}) \frac{\hat{\alpha}^2 \mathcal{N}(\bm{\uptheta};\bm{\uptheta}_i,\bm{\Lambda})}{|2\pi\bm{\Lambda}|^{-\frac{1}{2}}}  m(\bm{\uptheta}) d\bm{\uptheta} =  
    r \int m(\bm{\theta}) \mathcal{N}(\bm{\uptheta};\bm{\upmu}_r,\bm{\Sigma}_r) d\bm{\uptheta} = \\ \nonumber & r \bm{\upbeta}^{\sf T} \int \bm{\upphi}(\bm{\uptheta}) \mathcal{N}(\bm{\uptheta};\bm{\upmu}_r,\bm{\Sigma}_r) d\bm{\uptheta} = r \bm{\upbeta}^{\sf T} \mathbb{E}_{\bm{\uptheta}_r}\left[  \bm{\upphi}(\bm{\uptheta}_r) \right]
\end{align}
where $\mathbb{E}_{\bm{\uptheta}}\left[m(\theta) \mathbf{k}_{\zeta,\mathbf{z}}(\bm{\uptheta})\right] = \bm{\upmu}_{m \mathbf{k}_{\zeta,\mathbf{z}}}$, $r = \alpha^2 |\bm{\Lambda}^{-1} + \mathbf{I}|^{-\frac{1}{2}} \exp\left(-\frac{1}{2} \bm{\uptheta}_i^{\sf T} (\bm{\Lambda} + \mathbf{I})^{-1} \bm{\uptheta}_i \right)$, $\bm{\uptheta}_r \sim \mathcal{N}(\bm{\uptheta};\bm{\upmu}_r,\bm{\Sigma}_r)$ follows a multivariate Gaussian distribution, $\bm{\upmu}_r = \bm{\uptheta}_i(\bm{\Lambda} + \mathbf{I})$, and $\bm{\Sigma}_r = (\mathbf{I} + \bm{\Lambda}^{-1})^{-1}$. Note that $\bm{\upphi}(\bm{\uptheta}_r)$ is a vector of polynomial terms, for which the expectations are given by statistical moments, which have a closed form solution according to a multivariate Gaussian distribution see \citet{Wick1950}.    

\subsection*{Expectation of $\bm{\kappa}_{\zeta,\mathbf{z}}(\bm{\uptheta}) \bm{\kappa}_{\zeta,\mathbf{z}}^{\sf T}(\bm{\uptheta})$}
This term is again a cross-term and is dealt with in a similar way. 
\begin{align}
    & \mathbf{M}_{\bm{\kappa}_{\zeta,\mathbf{z}}\bm{\kappa}_{\zeta,\mathbf{z}}} =  \begin{bmatrix} 
\mathbb{E}_{\bm{\uptheta}}\left[\bm{\upphi}(\bm{\uptheta}) \bm{\upphi}^{\sf T}(\bm{\uptheta}) \right]  & \mathbb{E}_{\bm{\uptheta}}\left[ \bm{\upphi}(\bm{\uptheta}) \mathbf{k}^{\sf T}_{\zeta,\mathbf{z}}(\bm{\uptheta}) \right] \\
\mathbb{E}_{\bm{\uptheta}}\left[\mathbf{k}_{\zeta,\mathbf{z}}(\bm{\uptheta}) \bm{\upphi}^{\sf T}(\bm{\uptheta})  \right] & \mathbb{E}_{\bm{\uptheta}}\left[ \mathbf{k}_{\zeta,\mathbf{z}}(\bm{\uptheta}) \mathbf{k}_{\zeta,\mathbf{z}}^{\sf T}(\bm{\uptheta}) \right] 
\end{bmatrix}  
\end{align}
where $\mathbf{M}_{\bm{\kappa}_{\zeta,\mathbf{z}}\bm{\kappa}_{\zeta,\mathbf{z}}} = \mathbb{E}_{\bm{\uptheta}}\left[ \bm{\kappa}_{\zeta,\mathbf{z}}(\bm{\uptheta}) \bm{\kappa}_{\zeta,\mathbf{z}}^{\sf T}(\bm{\uptheta}) \right]$.

The bottom right expectation was determined previously. The top left can be expressed as follows:
\begin{align}
    \left[\mathbb{E}_{\bm{\uptheta}}\left[\bm{\upphi}(\bm{\uptheta}) \bm{\upphi}^{\sf T}(\bm{\uptheta}) \right]\right]_{ij} = \mathbb{E}_{\bm{\uptheta}} \left[\phi_i(\bm{\uptheta}) \phi_j(\bm{\uptheta}) \right] 
\end{align}
where the RHS expectation has a closed-form solution and is a common occurrence to determine the covariance employing PCE, see for example \citet{Dutta2010}.  

The last remaining term can be determined as follows:
\begin{align}
    & \nonumber \left[\mathbb{E}_{\bm{\uptheta}}\left[  \bm{\upphi}(\bm{\uptheta}) \mathbf{k}^{\sf T}_{\zeta,\mathbf{z}}(\bm{\uptheta}) \right] \right]_{ij} = \int \mathcal{N}(\bm{\uptheta};\mathbf{0},\mathbf{I}) \frac{\hat{\alpha}^2 \mathcal{N}(\bm{\uptheta};\bm{\uptheta}_i,\bm{\Lambda})}{|2\pi\bm{\Lambda}|^{-\frac{1}{2}}} \phi_j(\bm{\uptheta}) d\bm{\uptheta} = \\ & r \int \phi_j(\bm{\uptheta}) \mathcal{N}(\bm{\uptheta}_r;\bar{\bm{\upmu}}_r,\bm{\Sigma}_r) d\bm{\uptheta} = r \mathbb{E}_{\bm{\uptheta}_r}\left[\phi_j(\bm{\uptheta}_r)\right] 
\end{align}
where $\bm{\uptheta}_r \sim \mathcal{N}(\bm{\uptheta};\bm{\upmu}_r,\bm{\Sigma}_r)$ follows a multivariate Gaussian distribution. Note that $\phi_j(\bm{\uptheta}_r)$ is a multivariate polynomial, for which the expectations are given by statistical moments, which have a closed form solution according to a multivariate Gaussian distribution \citep{Wick1950}.
\end{document}